\documentclass[preprint,fleqn,sort&compress,3p]{elsarticle}
\usepackage{amsthm,amsmath}
\usepackage{amssymb}
\usepackage{mathrsfs}
\usepackage{amsfonts}

\usepackage{graphicx}
\usepackage{url}
\usepackage{bm}
\usepackage{latexsym}
\usepackage{hyperref}
\usepackage[hyphenbreaks]{breakurl}
\usepackage{booktabs}
\usepackage{txfonts}
\usepackage{color}

\newcommand{\cre}[1]{\color{red}{#1}}
\newtheorem{theorem}{Theorem}[section]
\newtheorem{proposition}{Proposition}[section]

\newtheorem{lemma}{Lemma}[section]

\numberwithin{equation}{section}
\newtheorem{remark}{Remark}[section]
\journal{Numerical Mathematics: Theory, Methods and Applications}

\begin{document}
\title{A note on parallel preconditioning for the all-at-once
solution of Riesz fractional diffusion equations}


\author[a,b]{Xian-Ming Gu}
\ead{guxianming@live.cn}
\author[c]{Yong-Liang Zhao\corref{cor1}}
\ead{ylzhaofde@sina.com}
\cortext[cor1]{Corresponding author}
\author[c]{Xi-Le Zhao}
\ead{xlzhao122003@163.com}
\author[d]{Bruno Carpentieri}
\ead{bcarpentieri@gmail.com}
\author[a]{Yu-Yun Huang}
\ead{18735189972@163.com}
\address[a]{School of Economic Mathematics/Institute of Mathematics, \\
Southwestern University of Finance and Economics, Chengdu, Sichuan 611130, P.R. China}
\address[b]{Bernoulli Institute for Mathematics, Computer Science and Artificial
Intelligence, \\University of Groningen, Nijenborgh 9, P.O. Box 407, 9700 AK Groningen, The Netherlands}
\address[c]{School of Mathematical Sciences,\\
University of Electronic Science and Technology of China, Chengdu, Sichuan 611731, P.R. China}
\address[d]{Facolt\`{a} di Scienze e Tecnologie informatiche,\\
Libera Universit\`{a} di Bolzano, Dominikanerplatz 3 - piazza Domenicani, 3 Italy
- 39100, Bozen-Bolzano.}



\begin{abstract}
The $p$-step backward difference formula (BDF) for solving systems of ODEs can be formulated as  all-at-once linear systems that are solved by  parallel-in-time preconditioned
Krylov subspace solvers (see McDonald, Pestana, and Wathen [SIAM J. Sci. Comput., 40(2) (2018):
A1012-A1033] and Lin and Ng [\href{https://arxiv.org/abs/2002.01108}{arXiv:2002.01108}, 2020, 17
pages]).
However, when the BDF$p$ ($2\leq p\leq6$) method is used to solve time-dependent PDEs, the generalization of these studies is not straightforward as $p$-step BDF is not selfstarting  for $p\geq 2$. In this note, we focus on the 2-step BDF
which is often superior to the trapezoidal rule for solving the Riesz fractional diffusion equations, and show that it results into an all-at-once discretized system that is a low-rank perturbation of a block triangular Toeplitz system. We first give an estimation of the condition number of the all-at-once systems and
then, capitalizing on previous work, we propose two block circulant (BC) preconditioners.
Both the invertibility of these two BC preconditioners and the eigenvalue distributions of
preconditioned matrices are discussed in details. An efficient implementation of
these BC preconditioners is also presented, including the fast computation
of dense structured Jacobi matrices. Finally, numerical experiments involving both the one- and two-dimensional
Riesz fractional diffusion equations are reported to support our theoretical findings.
\end{abstract}

\maketitle
\section{Introduction}
\label{sec1}
In this paper, we are particularly interested in the efficient numerical solution of evolutionary partial differential equations (PDEs) with both first order temporal derivative and space fractional-order derivative(s). These models arise
in various scientific applications in different fields including physics \cite{Bay2016}, bioengineering \cite{Magin06}, hydrology \cite{MMMCT}, and finance \cite{Scalas00}, etc., owing to the potential of fractional calculus to describe rather accurately natural processes which maintain long-memory and hereditary properties in complex systems~\cite{Kilbas,Magin06}. In particular, fractional diffusion equations can provide an adequate and accurate description of transport processes that exhibit anomalous diffusion, for example subdiffusive phenomena and L\'{e}vy fights \cite{Metz00}, which cannot be modelled properly by second-order diffusion equations. As most fractional diffusion equations can not be solved analytically, approximate numerical solutions are sought by using efficient numerical methods such as, e.g., (compact) finite difference~\cite{MMMCT,Celik,Chen2014,Tian15,Liao18,Hu2019}, finite element~\cite{Yue2019} and spectral (Galerkin) methods~\cite{Xu2019,Xu2020,Zhao20}.

Many numerical techniques proposed in the literature for solving this class of problems
are the common time-stepping schemes. They solve the underlying evolutionary PDEs with
space fractional derivative(s) by marching in time sequentially, one level after the other.
As many time steps may be usually necessary to balance the numerical errors arising from the spatial discretization, these conventional time-stepping schemes can be very time-consuming. This concern motivates the recent development of parallel-in-time (PinT) numerical solutions for evolutionary PDEs (especially with space fractional derivative(s)) including, e.g., the inverse Laplace transform method~\cite{Pang}, the MGRIT method~\cite{Gander16x,Falgout,Yue2019,Wu2019g}, the exponential integrator \cite{Farquhar} and the parareal method~\cite{Wu2017}. A class of PinT methods, i.e., the space-time methods,
solves the evolutionary PDEs at all the time levels simultaneously by performing an all-at-once discretization
that results into a large-scale linear system that is typically solved by preconditioned Krylov subspace methods; refer e.g., to \cite{Maday08,KYN,Podlub,Weaver,Gu2015,Lei17,Bertac,Breiten,Ran2019,Banjai12} for details. However, most
of them only focus on the numerical solutions of one-dimensional space fractional diffusion equations \cite{Gu2015,Lei17,Ran2019,Zhao19}
due to the huge computational cost required for high-dimensional problems.

Recently, McDonald, Pestana and Wathen proposed in~\cite{McDona} a block
circulant (BC) preconditioner to accelerate the convergence of Krylov subspace
methods for solving the all-at-once linear system arising from $p$-step BDF temporal
discretization of evolutionary PDEs. Parallel experiments with the BC preconditioner in~\cite{McDona} are reported by Goddard and Wathen in~\cite{Goddard}. In~\cite{Lin20}, a generalized version of the BC
preconditioner has been proposed by Lin and Ng who introduced a parameter $\alpha\in(0,1)$ into the top-right block of the BC preconditioner that can be  fine-tuned to handle efficiently the case of very small diffusion coefficients. Both the BC and the generalized BC preconditioners use a modified \emph{diagonalization technique} that is originally proposed by Maday and R{\o}nquist~\cite{Maday08,KYN}. The investigations in \cite{McDona,Goddard,Lin20,Wu2018} mainly focus on the 1-step BDF (i.e. the backward Euler method) for solving the underlying evolutionary PDEs, which results in an exact block lower triangular Toeplitz (BLTT) all-at-once system. On the other hand, when the BDF$p$ ($2\leq p\leq6$) method is used to discretize the evolutionary PDEs, the complete BLTT all-at-once systems cannot be obtained  as implicit schemes based on the BDF$p$ ($2\leq p\leq6$) for evolutionary PDEs are not selfstarting \cite{Xu2019,KBPS,Ascher,Bokano,Lubich13}. For example, when we establish the fully discretized scheme based on popular BDF2 for evolutionary PDEs, we often need to use the backward Euler method to compute the solution at the first time level \cite{Emmrich,Kucera,Nishik,Liao18}.

In this study, we particularly consider the second-order accurate implicit difference BDF2 scheme
for solving the one- and two-dimensional Riesz fractional diffusion equations (RFDEs).
Although the Crank-Nicolson (C-N) method \cite{Celik} is a very popular solution option for solving such
RFDEs, it is still only $A$-stable rather than $L$-stable \cite{Ascher,Zhang19}. By contrast,
the BDF2 scheme with stiff decay can be more ``stable" and slightly cheaper because it
is always unconditionally stable and the numerical solution is often guaranteed to be positive and physically more reliable near initial time for the numerical solutions of evolutionary PDEs with either integral or fractional order spatial derivatives; see \cite{Liao18,Bokano,Emmrich,Li2019x} for details. After the spatial discretization of Riesz fractional derivative(s), we reformulate the BDF2 scheme for
the semi-discretized system of ODEs into an all-at-once system, where its coefficient matrix is
a BLTT matrix with a low-rank perturbation. Then, we tightly estimate the condition number of the
all-at-once systems and adapt the generalized BC (also including the standard BC) preconditioner for such an all-at-once system. Meanwhile, the invertibility of the generalized
BC preconditioner is discussed apart from the work in \cite{McDona,Lin20},
and the eigenvalue distributions of preconditioned matrices dictating the convergence rate of Krylov subspace methods is preliminarily investigated.
From these
discussions, we derive clear arguments explaining the better performance of the generalized BC preconditioner against the BC
preconditioner, especially for very small diffusion coefficient.

The practical implementation of the generalized BC preconditioner requires to solve a sequence
of dense complex-shifted linear systems with real symmetric negative definite (block) Toeplitz
coefficient matrices, which arise from the numerical discretization of Riesz fractional derivatives.
Because solving such Toeplitz systems is often very time-consuming, classical circulant
preconditioners have been adapted to accelerate the convergence of iterative solutions, especially
for level-1 Toeplitz systems. However, circulant preconditioners are often less efficient to extend
for solving high-dimensional model problems (i.e., block Toeplitz discretized systems) \cite{Bertac,Nout11}.
In this work, by approximating the real symmetric (block) Toeplitz matrix by the (block) $\tau$-matrix~\cite{Bini90,Nout11}, that can be efficiently diagonalized using the discrete sine transforms (DSTs), we present an efficient implementation for the generalized BC preconditioner that does avoid any dense matrix storage and only needs to solve the first half of the
sequence of complex-shifted systems. Estimates of the computational cost and of the memory requirements of the associated generalized BC preconditioner are given. Our numerical experiments suggest that the BDF2 all-at-once system utilized the preconditioned Krylov subspace solvers can be a competitive solution method for RFDEs.

The rest of this paper is organized as follows. In the next section, the all-at-once linear system derived from the BDF2 scheme for solving the RFDEs is presented. Meanwhile, the invertibility
of the pertinent all-at-once system is proved and its condition number is estimated. In Section~\ref{sec3}, the generalized BC preconditioner is adapted and discussed, both the properties and the efficient implementation of the preconditioned system are analyzed. In Section \ref{sec4}, numerical results are reported to support our findings and the effectiveness of the proposed preconditioner. Finally, the paper closes with conclusions in Section~\ref{sec5}.

\section{The all-at-once system of Riesz fractional diffusion equations}
\label{sec2}
In this section, we present the development of a numerical scheme for initial-boundary
problem of Riesz fractional diffusion equations that preserves the positivity of the solutions. Then, in the next section, we discuss its efficient parallel implementation.
\subsection{The all-at-once discretization of Riesz fractional diffusion equation}
\label{sec2.1}
The governing Riesz fractional diffusion equations \footnote{Here it is worth noting
that if the Riesz fractional derivative is replaced with the fractional Laplacian,
we can just follow the work in \cite{Hao2019} for the spatial discretization of one- or
two-dimensional model problem and our proposed numerical techniques in the next sections are easy and simple to be adapted for such an extension.} of the anomalous diffusion
process can be written as
\begin{equation}
\begin{cases}
\frac{\partial u(x,t)}{\partial t} = \kappa_{\gamma}\frac{\partial^{\gamma}
u(x,t)}{\partial |x|^{\gamma}} + f(x,t),& (x,t)\in(a,b)\times(0,T],\\
u(x,0) = \phi(x),& x\in [a,b],\\
u(a,t) = u(b,t) = 0,& t\in[0,T],
\end{cases}
\label{eq1.1}
\end{equation}
where $u(x,t)$ may represent, for example, a solute concentration, and constant $\kappa_{\gamma} > 0$ the diffusion coefficient. This equations is a superdiffusive
model largely used in fluid flow analysis, financial modelling and other
applications. The Riesz fractional derivative $\partial^{\gamma}u(x,t)/\partial
|x|^{\gamma}$ is defined by \cite{Yang10}
\begin{equation}
\begin{split}
\frac{\partial^{\gamma}u(x,t)}{\partial |x|^{\gamma}} & =
-\frac{1}{2\cos\left(\frac{\pi\gamma}{2}\right)}\cdot\frac{1}{\Gamma(2 - \gamma)}\cdot
\frac{\partial^2}{\partial x^2}\int^{b}_a\frac{u(\xi,t)}{|x - \xi|^{\gamma-1}}d\xi\\
& = -\frac{1}{2\cos\left(\frac{\pi\gamma}{2}\right)}[
{}_aD^{\gamma}_xu(x,t) + {}_xD^{\gamma}_{b}u(x,t)],\quad \gamma\in(1,2),
\end{split}
\end{equation}
in which ${}_aD^{\gamma}_x$ and ${}_xD^{\gamma}_{b}$ are the left and right
Riemann-Liouville fractional derivatives of order $\gamma\in(1,2)$ given, respectively,
by \cite{Kilbas}
\begin{equation*}
{}_aD^{\gamma}_xu(x,t) = \frac{1}{\Gamma(2 - \gamma)}\frac{\partial^2}{\partial
x^2}\int^{x}_a\frac{u(s,t)}{(x - s)^{\gamma-1}}ds,
\end{equation*}
and
\begin{equation*}
{}_xD^{\gamma}_{b}u(x,t) = \frac{1}{\Gamma(2 - \gamma)}\frac{\partial^2}{\partial
x^2}\int^{b}_x\frac{u(s,t)}{(s - x)^{\gamma - 1}}ds.
\end{equation*}
As $\gamma\rightarrow 2$, it notes that Eq. (\ref{eq1.1})
degenerates into the classical diffusion equation \cite{Yang10}.

Next, we focus on the numerical solution of Eq. (\ref{eq1.1}). We consider
a rectangle $\bar{Q}_T = \{(x,t):~a\leq x\leq b,~0\leq t \leq T\}$ discretized
on the mesh $\varpi_{h\tau} = \varpi_h\times \varpi$, where $\varpi_h =\{x_i
= a + ih,~i=0,1,\cdots,N;~h=(b-a)/N\}$, and $\varpi_{\tau} = \{t_k = k\tau,~j
= 0,1,\cdots,M;~\tau=T/M\}$. We denote by $v =\{v_i~|~0\leq i\leq N\}$ any grid
function.

A considerable amount of work has been devoted in the past years to the development of fast methods for the approximation of the Riesz and Riemann-Liouville (R-L) fractional derivatives, such as the first-order
accurate shifted Gr{\"u}nwald approximation \cite{MMMCT,Yang10} and the second-order
accurate weighted-shifted Gr{\"u}nwald difference (WSGD) approximation \cite{Tian15,Hao2015}.
Due to the relationship between the two kinds of fractional derivatives, all the numerical schemes proposed for the R-L fractional derivatives can be easily adapted to approximate the Riesz
fractional derivative. Although the solution approach described in this work can accomodate
any spatial discretized method, we choose the so-called fractional centred difference
formula \cite{Celik} of the Riesz fractional derivatives for clarity.

For any function $u(x)\in L^1(\mathbb{R})$, we denote
\begin{equation}
\Delta^{\gamma}_hu(x) = -\frac{1}{h^{\gamma}}\sum^{[(x-a)/h]}_{\ell = -[(b-x)/h]}
\omega^{(\gamma)}_{\ell}u(x - \ell h),\quad x\in\mathbb{R},
\label{eq2.3}
\end{equation}
where the $\gamma$-dependent weight coefficient is defined as
\begin{equation}
\omega^{(\gamma)}_{\ell} = \frac{(-1)^{\ell}\Gamma(1 + \gamma)}{
\Gamma(1 + \gamma/2 - \ell)\Gamma(1 + \gamma/2 + \ell)},\quad \ell
\in\mathbb{Z}.
\label{eq2.4}
\end{equation}
As noted in \cite{Celik}, $\omega^{(\gamma)}_{\ell} = \mathcal{O}(\ell^{-1-
\gamma})$ and the fractional centred difference formula $\Delta^{\gamma}_hu(x)$
exists for any $u(x)\in L^1(\mathbb{R})$. Some properties of the coefficient
$\omega^{(\gamma)}_{\ell}$ and the operator $\Delta^{\gamma}_hu(x)$ are presented in the following lemmas.
\begin{lemma}{\rm (\cite{Celik})}
For $\gamma\in(1,2)$, the coefficient $\omega^{(\gamma)}_{\ell}$, $\ell\in\mathbb{Z}$,
defined in (\ref{eq2.4}) fulfils
\begin{equation}
\begin{cases}
\omega^{(\gamma)}_0\geq 0,\quad \omega^{(\gamma)}_{-\ell} =
\omega^{(\gamma)}_{\ell},\quad |\ell|\geq 1,\\
\sum\limits^{\infty}_{\ell = -\infty}\omega^{(\gamma)}_{\ell} = 0,\quad
\omega^{(\gamma)}_0 = \sum\limits^{-1}_{\ell = -\infty}
|\omega^{(\gamma)}_{\ell}| + \sum\limits^{\infty}_{\ell = 1}
|\omega^{(\gamma)}_{\ell}|.
\end{cases}
\end{equation}
\end{lemma}
\begin{lemma}
{\rm(\cite{Celik,Hao2019})} Suppose that $u\in L^1(\mathbb{R})$ and
\begin{equation*}
  u(x) \in \mathscr{C}^{2+\gamma}(\mathbb{R}):=\bigg\{u\bigg|\int\nolimits_{
  -\infty}^{+\infty}(1+|\xi|)^{2+\gamma}|\hat{u}(\xi)|d\xi<\infty\bigg\},
\end{equation*}
where $\hat{u}(\xi)$ is the Fourier transformation of $u(x)$. Then for a fixed
$h$, the fractional centred difference operator in (\ref{eq2.3}) holds
\begin{equation*}
\frac{\partial^{\gamma}u(x)}{\partial |x|^{\gamma}} = \Delta^{\gamma}_hu(x) +
\mathcal{O}(h^2)
\end{equation*}
uniformly for $x\in \mathbb{R}$ and $u(x)\equiv0$ ($x\in\mathbb{R}\backslash[a,b]$).
In particular, if $\gamma=2$, then it coincides with the
second-order derivative  approximation.
\label{lem2.2}
\end{lemma}

At this stage, let $u(x,t)\in\mathcal{C}^{4,3}_{x,t}([a,b]\times[0,T])$ be a solution to the problem (\ref{eq1.1}) and consider Eq. (\ref{eq1.1}) at the set of grid points
$(x,t)=(x_i,t_k)\in\bar{Q}_T$ with $i=1,2,\cdots,N-1,~k=1,\cdots,M$. The first-order
time derivative at the point $t = t_k$ is approximated by the second-order backward difference formula, i.e.,
\begin{equation}
\frac{\partial u(x,t)}{\partial t}\Big|_{t = t_k} = \frac{u(x,t_k) - 4u(x,t_{k-1})
+3u(x,t_{k-2})}{2\tau} + \mathcal{O}(\tau^2),\quad k \geq 2,
\end{equation}
then we define $U^{k}_i = u(x_i,t_k)$ and $f^{k}_i = f(x_i,t_k)$ to obtain
\begin{equation}
\begin{cases}
\frac{U^{k}_i - 4U^{k-1}_i + 3U^{k-2}_i}{2\tau} = \kappa_{\gamma}\Delta^{
\gamma}_hU^{k}_i + f^{k}_i + R^{k}_i,&1\leq i\leq N-1,~2\leq k\leq M,\\
U^{0}_i = \phi(x_i), & 1\leq i\leq N-1,\\
U^{k}_0 = U^{k}_N = 0,& 0\leq k\leq M,
\end{cases}
\label{eq2.7}
\end{equation}
where $\{R^{k}_i\}$ are small and satisfy the inequality
\begin{equation*}
|R^{k}_i|\leq c(\tau^2 + h^2),\quad 1\leq i\leq N-1,~2\leq k\leq M.
\end{equation*}

Note that we cannot omit the small term in the derivation of a two-level difference scheme that is not selfstarting from Eq.~(\ref{eq2.7}) due to the unknown information of $u(x_i,
t_1)$. One of the most popular strategies to compute the first time step solution $u^{1}_i$ is to use a first-order backward Euler scheme (which is only used once). This yields the following implicit difference scheme \cite{Bokano,Liao18} for Eq. (\ref{eq1.1}):
\begin{equation}
\begin{cases}
D^{2}_tu^{k}_i = \kappa_{\gamma}\Delta^{\gamma}_hu^{k}_i + f^{k}_i,&1\leq i\leq N-1,
~2\leq k\leq M,\\
u^{0}_i = \phi(x_i), & 1\leq i\leq N-1,\\
u^{k}_0 = u^{k}_N = 0,& 0\leq k\leq M,
\end{cases}
\end{equation}
where
\begin{equation*}
D^{2}_tu^{k}_i =
\begin{cases}
\frac{u^{k}_i - 4u^{k-1}_i + 3u^{k-2}_i}{2\tau},& 2\leq k\leq M,\\
\frac{u^{1}_i - u^{0}_i}{\tau}, & k = 1.
\end{cases}
\end{equation*}

In order to implement the proposed scheme, here we define ${\bm u}^{k} = [u(x_1,t_k),u(x_2,t_k),
\cdots,u(x_{N-1},t_k)]^\top$ and ${\bm f}^{k} = [f(x_1,t_k),f(x_2,t_k),\cdots,f(x_{N-1},
t_k)]^\top$; then $\kappa_{\gamma}\Delta^{\gamma}_hu^{k}_i$ can be written into the
matrix-vector product form $A{\bm u}^{k}$, where
\begin{equation}
A = -\kappa_{\gamma}T_x =
-\frac{\kappa_{\gamma}}{h^\gamma}
\begin{bmatrix}
\omega^{(\gamma)}_0 &\omega^{(\gamma)}_{-1} &\omega^{(\gamma)}_{-2} &\cdots
&\omega^{(\gamma)}_{3-N}&\omega^{(\gamma)}_{2-N}\\
\omega^{(\gamma)}_1 &\omega^{(\gamma)}_0&\omega^{(\gamma)}_{-1}&\cdots
&\omega^{(\gamma)}_{4-N}&\omega^{(\gamma)}_{3-N}\\
\omega^{(\gamma)}_2 &\omega^{(\gamma)}_1&\omega^{(\gamma)}_{0}&\cdots
&\omega^{(\gamma)}_{5-N}&\omega^{(\gamma)}_{4-N}\\
\vdots&\vdots&\vdots&\ddots&\vdots&\vdots\\
\omega^{(\gamma)}_{N-3} &\omega^{(\gamma)}_{N-4}&\omega^{(\gamma)}_{N-5}&\cdots
&\omega^{(\gamma)}_0&\omega^{(\gamma)}_{-1}\\
\omega^{(\gamma)}_{N-2} &\omega^{(\gamma)}_{N-3}&\omega^{(\gamma)}_{N-4}&\cdots
&\omega^{(\gamma)}_1&\omega^{(\gamma)}_0
\end{bmatrix}\in\mathbb{R}^{(N-1)\times(N-1)}.
\end{equation}
It is easy to prove that $T_x$ is a symmetric positive definite (SPD) Toeplitz matrix
(see \cite{Celik}). Therefore, it can be stored with only $N - 1$ entries and the fast
Fourier transforms (FFTs) can be applied to carry out the matrix-vector product in
$\mathcal{O}(N\log N)$ operations.
The matrix-vector form of the 2-step BDF method with start-up backward Euler method for solving the model problem (\ref{eq1.1}) can
be formulated as follows:
\begin{subequations}
\begin{align}
\frac{{\bm u}^{1} - {\bm u}^{0}}{\tau} - A{\bm u}^{1}
& = {\bm f}^{1},\label{eq1.2a} \\
\frac{3{\bm u}^{k} - 4{\bm u}^{k-1} + {\bm u}^{k-2}}{2\tau} -
A{\bm u}^{k} & = {\bm f}^{k},\quad 2\leq k\leq N_t.
\label{eq1.2b}
\end{align}
\label{eq1.2}
\end{subequations}
We refer to matrix $A$ as the Jacobian matrix. In particular, we note that the numerical scheme (\ref{eq1.2}) has been proved to be unconditionally stable and second-order convergent in the discrete $L^2$ norm~\cite{Liao18}.
Note that it is computationally more efficient than the C-N scheme presented in~\cite{Celik}
as it requires one less matrix-vector multiplication per time step.

Instead of computing the solution of (\ref{eq1.2}) step-by-step, we try to get an all-at-once approximation by solving the following linear system:
\begin{equation}
\left(\begin{bmatrix}
1  &  \\
-2 & \frac{3}{2} \\
\frac{1}{2} &-2 & \frac{3}{2} \\
   &\ddots & \ddots & \ddots \\
   &        &\frac{1}{2}& -2         & \frac{3}{2} & \\
   &        &            &\frac{1}{2}&-2           &\frac{3}{2}
\end{bmatrix}\otimes I_s - \tau I_t\otimes A\right)
\begin{bmatrix}
{\bm u}^1 \\
{\bm u}^2 \\
\vdots \\
{\bm u}^{N_t - 1} \\
{\bm u}^{N_t} \\
\end{bmatrix} =
\begin{bmatrix}
\tau{\bm f}^1 + {\bm u}^{0}\\
\tau{\bm f}^2 - \frac{{\bm u}^{0}}{2}\\
\vdots \\
\tau{\bm f}^{N_t - 1} \\
\tau{\bm f}^{N_t} \\
\end{bmatrix},
\label{eq11}
\end{equation}
where $I_s$ and $I_t$ are two identity matrices of order $N_s~(=N-1)$ and $N_t~(=M)$, respectively.
We denote the above linear system as follows:
\begin{equation}
\mathcal{A}{\bm U} = {\bm F},\quad {\bm U} = \begin{bmatrix}
{\bm u}^1 \\
{\bm u}^2 \\
\vdots \\
{\bm u}^{N_t - 1} \\
{\bm u}^{N_t} \\
\end{bmatrix},\quad {\bm F} =
\begin{bmatrix}
\tau{\bm f}^1 + {\bm u}^{0}\\
\tau{\bm f}^2 - \frac{{\bm u}^{0}}{2}\\
\vdots \\
\tau{\bm f}^{N_t - 1} \\
\tau{\bm f}^{N_t} \\
\end{bmatrix},
\label{eq12}
\end{equation}
where $\mathcal{A} = C\otimes I_s -\tau I_t\otimes A$ with
\begin{equation}
C = \begin{bmatrix}
1  &  \\
-2 & \frac{3}{2} \\
\frac{1}{2} &-2 & \frac{3}{2} \\
   &\ddots & \ddots & \ddots \\
   &        &\frac{1}{2}& -2         & \frac{3}{2} & \\
   &        &            &\frac{1}{2}&-2           &\frac{3}{2}
\end{bmatrix}\in\mathbb{R}^{N_t\times N_t},
\label{eq1.7}
\end{equation}
and it is clear that $\mathcal{A}$ is invertible, because its all diagonal blocks
(i.e, either $I_s - \tau A$ or $\frac{3}{2}I_s - \tau A$) are invertible \cite{McDona,Goddard}.

\subsection{Properties of the all-at-once system}
In this subsection, we investigate the properties of the discrete all-at-once
formulation (\ref{eq11}). This will guide us to discuss the design of an efficient solver for such
a large linear system.

\begin{lemma}
For the matrix $C$ in (\ref{eq1.7}), we have the following estimates,
\begin{itemize}
  \item[1)] $\|C^{-1}\|_{\infty} \leq \frac{3N_t}{2}$;
  \item[2)] $\|C^{-1}\|_1 = N_t$.
\end{itemize}
\label{lem2.3}
\end{lemma}
\textbf{Proof}. Consider the following matrix splitting and then perform the matrix factorization,
\begin{equation}
\begin{split}
C  & =
\begin{bmatrix}
\frac{3}{2}  &  \\
-2 & \frac{3}{2} \\
\frac{1}{2} &-2 & \frac{3}{2} \\
   &\ddots & \ddots & \ddots \\
   &        &\frac{1}{2}& -2         & \frac{3}{2} & \\
   &        &            &\frac{1}{2}&-2           &\frac{3}{2}
\end{bmatrix}
-
\begin{bmatrix}
\frac{1}{2}  &  \\
0 & 0 \\
0 &0 & 0 \\
   &\ddots & \ddots & \ddots \\
   &        &0& 0         & 0 & \\
   &        &            &0&0           &0
\end{bmatrix} \triangleq \hat{C} - \frac{1}{2}{\bm e}_1{\bm e}^{T}_1\\
&= \frac{3}{2}
\begin{bmatrix}
1&   \\
-\frac{1}{3}&1 \\
&-\frac{1}{3}&1 \\
& &\ddots&\ddots\\
& & &-\frac{1}{3}&1 & \\
& & & & -\frac{1}{3} &1
\end{bmatrix}
\cdot
\begin{bmatrix}
1&   \\
-1&1 \\
&-1&1 \\
& &\ddots&\ddots\\
& & &-1&1 & \\
& & & & -1 &1
\end{bmatrix} - \frac{1}{2}{\bm e}_1{\bm e}^{T}_1
\end{split}
\label{eq1.9}
\end{equation}
where the vector ${\bm e}_1 = [1,0,\cdots,0]^T\in\mathbb{R}^{N_t}$. According to
the Sherman-Morrison formula \cite{Hager89}, we can write
\begin{equation}
C^{-1} = \hat{C}^{-1} + \frac{\hat{C}^{-1}\frac{1}{2}{\bm e}_1{\bm e}^{T}_1
\hat{C}^{-1}}{1 - \frac{1}{2}{\bm e}^{T}_1\hat{C}^{-1}{\bm e}_1} =
\hat{C}^{-1} + \frac{\hat{C}^{-1}{\bm e}_1{\bm e}^{T}_1\hat{C}^{-1}}{2 -
{\bm e}^{T}_1\hat{C}^{-1}{\bm e}_1}.
\end{equation}
On the other hand, we can compute the inverse of $\hat{C}$ as follows,
\begin{equation}
\hat{C}^{-1}
= \frac{2}{3}
\begin{bmatrix}
1      &                  \\
1      & 1     &          \\
1      & 1     & 1     &  \\
\vdots & \ddots &\ddots &\ddots &  \\
1      &1  &\cdots &1      & 1  & \\
1      &1  &\cdots &1      & 1  & 1
\end{bmatrix}\cdot
\begin{bmatrix}
1&\\
\frac{1}{3}           & 1                \\
\frac{1}{3^2}         & \frac{1}{3} & 1  \\
\vdots                &\ddots       &\ddots  &\ddots \\
\frac{1}{3^{N_t - 2}} &\frac{1}{3^{N_t - 3}} &\cdots &\frac{1}{3}   & 1         &\\
\frac{1}{3^{N_t - 1}}     &\frac{1}{3^{N_t - 2}} &\cdots &\frac{1}{3^2} &\frac{1}{3}&1  \\
\end{bmatrix}.
\end{equation}
Therefore, we know that
\begin{equation}
\hat{C}^{-1}{\bm e}_1 =
\frac{2}{3}
\begin{bmatrix}
1      &                  \\
1      & 1     &          \\
1      & 1     & 1     &  \\
\vdots & \ddots &\ddots &\ddots &  \\
1      &1  &\cdots &1      & 1  & \\
1      &1  &\cdots &1      & 1  & 1
\end{bmatrix}\cdot
\begin{bmatrix}
1&\\
\frac{1}{3}           \\
\frac{1}{3^2}         \\
\vdots                \\
\frac{1}{3^{N_t - 2}} \\
\frac{1}{3^{N_t - 1}}     \\
\end{bmatrix}
=
\begin{bmatrix}
\frac{2}{3}&\\
1 - \frac{1}{3^2} \\
1 - \frac{1}{3^3} \\
\vdots                \\
1 - \frac{1}{3^{N_t - 1}} \\
1 - \frac{1}{3^{N_t}}     \\
\end{bmatrix}
\end{equation}
and
\begin{equation}
{\bm e}^{T}_1\hat{C}^{-1} = \frac{2}{3}{\bm e}^{T}_1,\quad
{\bm e}^{T}_1\hat{C}^{-1}{\bm e}_1 = \frac{2}{3}.
\end{equation}
The explicit expression of $C^{-1}$ has the following form:
\begin{equation}
C^{-1} = \hat{C}^{-1} + \frac{1}{2}
\begin{bmatrix}
\frac{2}{3}&\\
1 - \frac{1}{3^2} \\
1 - \frac{1}{3^3} \\
\vdots                \\
1 - \frac{1}{3^{N_t - 1}} \\
1 - \frac{1}{3^{N_t}}     \\
\end{bmatrix}{\bm e}^{T}_1
=
\hat{C}^{-1} + \frac{1}{2}
\begin{bmatrix}
\frac{2}{3}               &0      &\cdots &0\\
1 - \frac{1}{3^2}         &0      &\cdots &0\\
1 - \frac{1}{3^3}         &0      &\cdots &0\\
\vdots                    &\vdots &\cdots &\vdots\\
1 - \frac{1}{3^{N_t - 1}} &0      &\cdots &0\\
1 - \frac{1}{3^{N_t}}     &0      &\cdots &0\\
\end{bmatrix}.
\end{equation}
Since $\hat{C}^{-1}$ is a lower triangular Toeplitz matrix, $\|C^{-1}\|_1$ is the absolute sum of elements of its last row, i.e.,
\begin{equation}
\begin{split}
\|C^{-1}\|_1 & = \frac{2}{3}\left[\left(1 + \cdots + \frac{1}{3^{N_t-1}}\right) +
\left(1 + \cdots + \frac{1}{3^{N_t-2}}\right) + \cdots + \left(1 + \frac{1}{3}\right)
+1 \right] + \frac{1}{2}\left(1 - \frac{1}{3^{N_t}}\right)\\
& = N_t - \left(\frac{1}{3} + \frac{1}{3^2} + \cdots + \frac{1}{3^{N_t}}\right) +
\frac{1}{2}\left(1 - \frac{1}{3^{N_t}}\right)\\
& = N_t.
\end{split}
\end{equation}
Analogously, $\|C^{-1}\|_{\infty}$ is the absolute sum of elements of its
first column, i.e.,
\begin{equation}
\begin{split}
\|C^{-1}\|_{\infty} & = \frac{2}{3}\left[1 + \left(1 + \frac{1}{3}\right) + \cdots +
\left(1 + \cdots + \frac{1}{3^{N_t-2}}\right) + \left(1 + \cdots + \frac{1}{3^{N_t-1}}\right)
\right]\\
& + \frac{1}{2}\left(1 - \frac{1}{3} + 1 - \frac{1}{3^2} +
\cdots + 1 - \frac{1}{3^{N_t}}\right)\\
& = \frac{2}{3}\left(\frac{1 - \frac{1}{3^1}}{1 - \frac{1}{3}} + \frac{1 - \frac{1}{3^2}}{1
- \frac{1}{3}} + \cdots + \frac{1 - \frac{1}{3^{N_t-1}}}{1 - \frac{1}{3}} + \frac{1 - \frac{1}
{3^{N_t}}}{1 - \frac{1}{3}} + \right) + \frac{1}{2}\left(N_t - \frac{1 - \frac{1}{3^{N_t}
}}{1 - \frac{1}{3}}\cdot\frac{1}{3}\right)\\
& = \frac{3N_t}{2} - \frac{3}{4}\left(1 - \frac{1}{3^{N_t}}\right)\\
& \leq \frac{3N_t - 1}{2} \leq \frac{3N_t}{2}.
\end{split}
\end{equation}
Therefore, the above estimates are proved. \hfill$\Box$
\begin{theorem}~\label{theo:cond}
Suppose that $Re(\lambda(A))\leq 0$. Then, the following bounds hold for the critical singular values of
$\mathcal{A}$ in (\ref{eq11}):
\begin{equation}
\sigma_{max}(\mathcal{A})\leq 4 + \tau\|A\|_2\quad{\rm and}\quad
\sigma_{min}(\mathcal{A})\geq \frac{\sqrt{6}}{3N_t},
\end{equation}
so that ${\rm cond}(\mathcal{A})\leq 2\sqrt{6}N_t + \frac{\sqrt{6}T\|A\|_2}{2}$.
\label{thm2.1}
\end{theorem}
\textbf{Proof}. Since $Re(\lambda(A))\leq 0$, we may claim that $\|\mathcal{A}
{\bm z}\|\geq\|(C\otimes I_s){\bm z}\|$ for any vector ${\bm z}$. In particular, using the properties of the Kronecker product, we may estimate $\sigma_{min
}(\mathcal{A})\geq \sigma_{min}(C)$. The latter is computed using the spectral
norm of $C^{-1}$, which can be bounded by the following inequality:
\begin{equation}
\|C^{-1}\|_2\leq \sqrt{\|C^{-1}\|_1\|C^{-1}\|_{\infty}}.
\end{equation}
According to Lemma \ref{lem2.3}, we have
\begin{equation}
\sigma_{min}(\mathcal{A})\geq \sigma_{max}(C) =
\frac{1}{\|C^{-1}\|_2}\geq \frac{1}{\sqrt{\|C^{-1}\|_1\|C^{-1}\|_{\infty}}}
\geq \frac{\sqrt{6}}{3 N_t},
\end{equation}
which proves the second bound.

The first estimate is derived similarly. By applying the triangle inequality, we can write
\begin{equation}
\sigma_{max}(\mathcal{A}) = \|\mathcal{A}\|_2 \leq
\|C\|_2 + \tau\|A\|_2,
\end{equation}
whereas
\begin{equation}
\|C\|_1 = \|C\|_{\infty} = 4,
\end{equation}
so the spectral norms are not greater than 4. Finally, the proof is completed by recalling
that $\tau = T/N_t$. \hfill$\Box$

The condition number can be expected to vary linearly with the main properties of the system, such as number of time steps, length of time interval, and norm of the Jacobian matrix~\cite{KBPS,Ascher,Gu2020}. More refined estimates could be provided by taking into account particular properties of $A$. In Theorem~\ref{theo:cond},  the time interval $[0,T]$   may not necessarily be equal to the whole observation range of an application. The global time scheme~(\ref{eq11}) could be applied by splitting the desired interval $[0, \hat{T}]$ into a sequence of subintervals $[0, T], [T, 2T],\cdots , [\hat{T}-T, \hat{T}]$ for ``large" time steps of size $T$ each,  solving by~(\ref{eq11}) for $[(q - 1)T, qT],~q = 1,\cdots, \hat{T}/T$, extracting the last snapshot $x_{Nt}$ and \textit{restarting} the method using $x_{Nt}$ as the initial state for the next interval. The optimal value of $T$ should provide the fastest computation~\cite{KBPS}. Another way to accelerate the solution is to reduce the condition number of the linear system by \textit{preconditioning}. This standard computational aspect will be considered in the next section.


%

\section{Parallel-in-time (PinT) preconditioners}
\label{sec3}
According to Theorem \ref{thm2.1}, when an iterative method, namely a Krylov subspace method,
is used for solving the all-at-once system (\ref{eq11}), it can converge very slowly or even stagnate. Therefore, in this section we look at preconditioners that can be efficiently implemented in the PinT framework.

\subsection{The structuring-circulant preconditioners}
\label{sec3.1}
Since the matrix $C$ defined in~(\ref{eq1.7}) can be viewed as the sum of a Toeplitz matrix and a rank-1 matrix, it is natural to define our first structuring-circulant preconditioner as
\begin{equation}
\mathcal{P}_{\alpha} = C_{\alpha}\otimes I_s - \tau I_t\otimes A
\end{equation}
where
\begin{equation}
C_{\alpha} =
\begin{bmatrix}
\frac{3}{2}  &  &&&\frac{\alpha}{2}&-2\alpha\\
-2 & \frac{3}{2} &&&&\frac{\alpha}{2}\\
\frac{1}{2} &-2 & \frac{3}{2} \\
   &\ddots & \ddots & \ddots \\
   &        &\frac{1}{2}& -2         & \frac{3}{2} & \\
   &        &            &\frac{1}{2}&-2           &\frac{3}{2}
\end{bmatrix},\quad \alpha\in(0,1],
\label{eq1.11}
\end{equation}
is a $\alpha$-circulant matrix of Strang type that can be diagonalized as
\begin{equation}
C_\alpha=V_\alpha D_\alpha V^{-1}_\alpha
\end{equation}
with
  \begin{subequations}
\begin{equation}
V_\alpha = \Lambda_\alpha\mathbb{F}_{N_t}^*,\quad~
D_\alpha = \text{diag}\big(\sqrt{N_t}\mathbb{F}_{N_t}
\Lambda_{\alpha}^{-1}C_\alpha(:,1)\big)
 = {\rm diag}\big(\lambda^{(\alpha)}_1,\lambda^{(\alpha)}_2,
\cdots,\lambda^{(\alpha)}_{N_t}\big)
 \label{eq1.6a}
 \end{equation}
where $\lambda^{(\alpha)}_n = \sum^{2}_{j = 0}r_j\alpha^{j/N_t}\theta^{(n-1)j}
~(r_0 = 3/2,~r_1 = -2,~r_2 = 1/2)$, `$*$' denotes the conjugate transpose of a
matrix, $C_\alpha(:,1)$ is the first column of $C_\alpha$, $\Lambda_{\alpha} =
{\rm diag}\big(1,\alpha^{-\frac{1}{N_t}},\cdots,\alpha^{-\frac{N_t-1}{N_t}}\big)
$ and
\begin{equation}
\begin{split}
\mathbb{F}_{N_t}=\frac{1}{\sqrt{N_t}} \begin{bmatrix}
 1 &1 &\dots &1\\
 1 &\theta &\dots &\theta^{N_t-1}\\
 \vdots &\vdots &\dots &\vdots\\
 1 &\theta^{N_t-1} &\dots &\theta^{(N_t-1)(N_t-1)}
  \end{bmatrix},\quad~ \theta=e^{\frac{2\pi {\rm i}}{N_t}}.
  \end{split}
  \label{eq1.6b}
 \end{equation}
\end{subequations}
The matrix $\mathbb{F}_{N_t}$ in~(\ref{eq1.6b}) is the discrete Fourier matrix,
and is unitary. In addition, it is worth noting that if we choose
$\alpha$ very small, then the condition number of $V_{\alpha}$ increases quickly
and therefore large roundoff errors may be introduced in our proposed method; see \cite{Wu2018,Wu2019g}
for a discussion of these issues. Using the property of the Kronecker product, we can factorize $\mathcal{P}_{\alpha}$  as
\begin{equation}
\mathcal{P}_{\alpha} = (V_\alpha\otimes I_s)\left(D_\alpha\otimes I_s - \tau
I_t\otimes A\right)(V^{-1}_{\alpha}\otimes I_s),
\label{eq1.7x}
\end{equation}
and this implies that we can compute ${\bm z} = \mathcal{P}^{-1}_{\alpha}{\bm v}$ via the following three steps:
\begin{equation*}
\begin{cases}
{\bm z}_1 = (V^{-1}_\alpha\otimes I_s){\bm v}
=(\mathbb{F}_{N_t}\otimes I_s)\left[(\Lambda^{-1}_{\alpha}\otimes I_s){\bm v}\right],
& \text{Step-(a)}, \\
\big(\lambda^{(\alpha)}_{n} I_s - \tau A\big){\bm z}_{2,n}= {\bm z}_{1,n},\quad  ~n=1,2,\dots,  {N_t},  & \text{Step-(b)}, \\
{\bm z} = (V_\alpha\otimes I_s){\bm z}_2=(\Lambda_{\alpha}\otimes I_s)\left[(\mathbb{F}^*_{N_t}\otimes I_s){\bm z}_{2}\right],  & \text{Step-(c)}, \\
 \end{cases}
\label{eq1.8}
 \end{equation*}
where ${\bm z}_{j}=({\bm z}_{j,1}^\top, {\bm z}_{j,2}^\top,\dots,{\bm z}_{j,N_t}^\top)^\top$ (with $j=1,2$) and $\lambda^{(\alpha)}_n$ is the $n$-th eigenvalue of $C_\alpha$. The first and third steps only involve matrix-vector multiplications with $N_t$ subvectors, and these can be computed simultaneously by the FFTs\footnote{It notes that there are also parallel FFTs available at \url{http://www.fftw.org/parallel/parallel-fftw.html}.} in $N_t$ CPUs \cite{Gu2020}.
The major computational cost is spent in the second step, however this step is naturally parallel for all the $N_t$ time steps.

Next, we study the invertibility of the matrix $\mathcal{P}_{\alpha}$. According to Eq. (\ref{eq1.7})
and because matrix $A$ is negative definite, the following result assures the invertibility of
$\mathcal{P}_{\alpha}$\footnote{It is unnecessary to keep that $r_0 \geq |r_1| + |r_2|$, which is introduced in \cite[Remark 4]{Lin20}.}.
\begin{proposition}
For $\alpha\in(0,1]$ and $\lambda^{(\alpha)}_n\in\mathbb{C}$, it holds that $\mathcal{R}e(\lambda^{(\alpha)}_n)\geq 0$,
where the equality is true if and only if $\alpha = 1$.
\label{pro3.1}
\vspace{-4mm}
\end{proposition}
\noindent\textbf{Proof}. Set $\varepsilon = \alpha^{1/N_t}\theta^{(n-1)} \triangleq
x + {\rm i}y$, then $x^2 + y^2 = \alpha^{2/N_t}\leq 1$, we have
\begin{equation}
\begin{split}
\mathcal{R}e(\lambda^{(\alpha)}_n) & = \mathcal{R}e\left(\frac{3}{2} - 2\varepsilon + \frac{\varepsilon^2}{2}\right) \\
& = \mathcal{R}e\left(\frac{1}{2}(1 - x - {\rm i}y) (3 - x - {\rm i}y)\right)\\
& = \frac{1}{2}[(1 - x)(3 - x) - y^2]\\
& \geq x^2 - 2x + 1 \\
& \geq 0,
\end{split}
\end{equation}
where the penultimate inequality holds due to $x^2 - 1\leq -y^2$. Moreover, if
$\alpha\neq1$, then it should hold that $x^2 + y^2 < 1$, i.e., $-1< x < 1$ and $x^2 -2x +
1 > 0$, which completes the proof. \hfill$\Box$

\begin{figure}[!htpb]
\centering
\includegraphics[width=2.137in,height=1.9in]{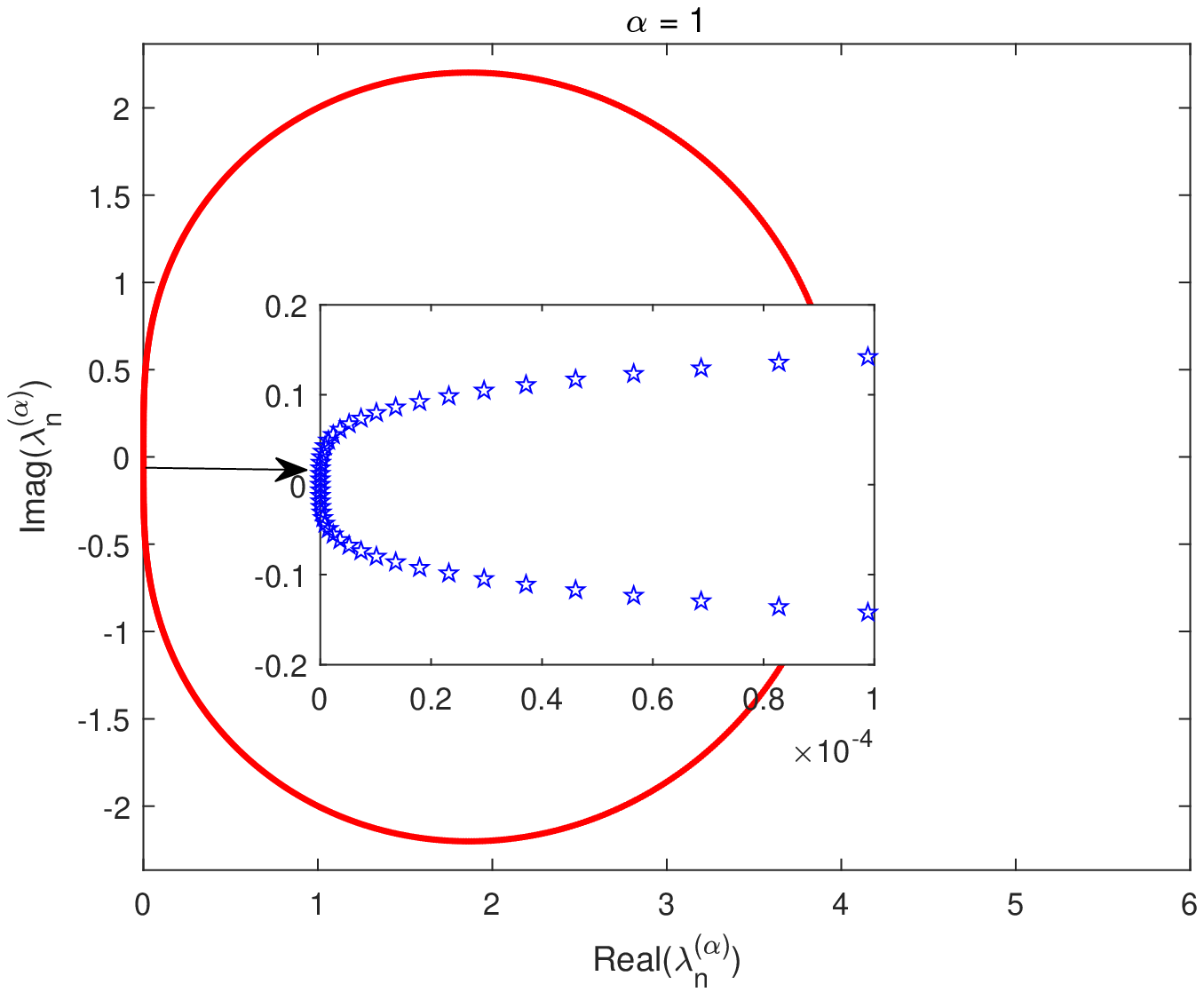}
\includegraphics[width=2.137in,height=1.9in]{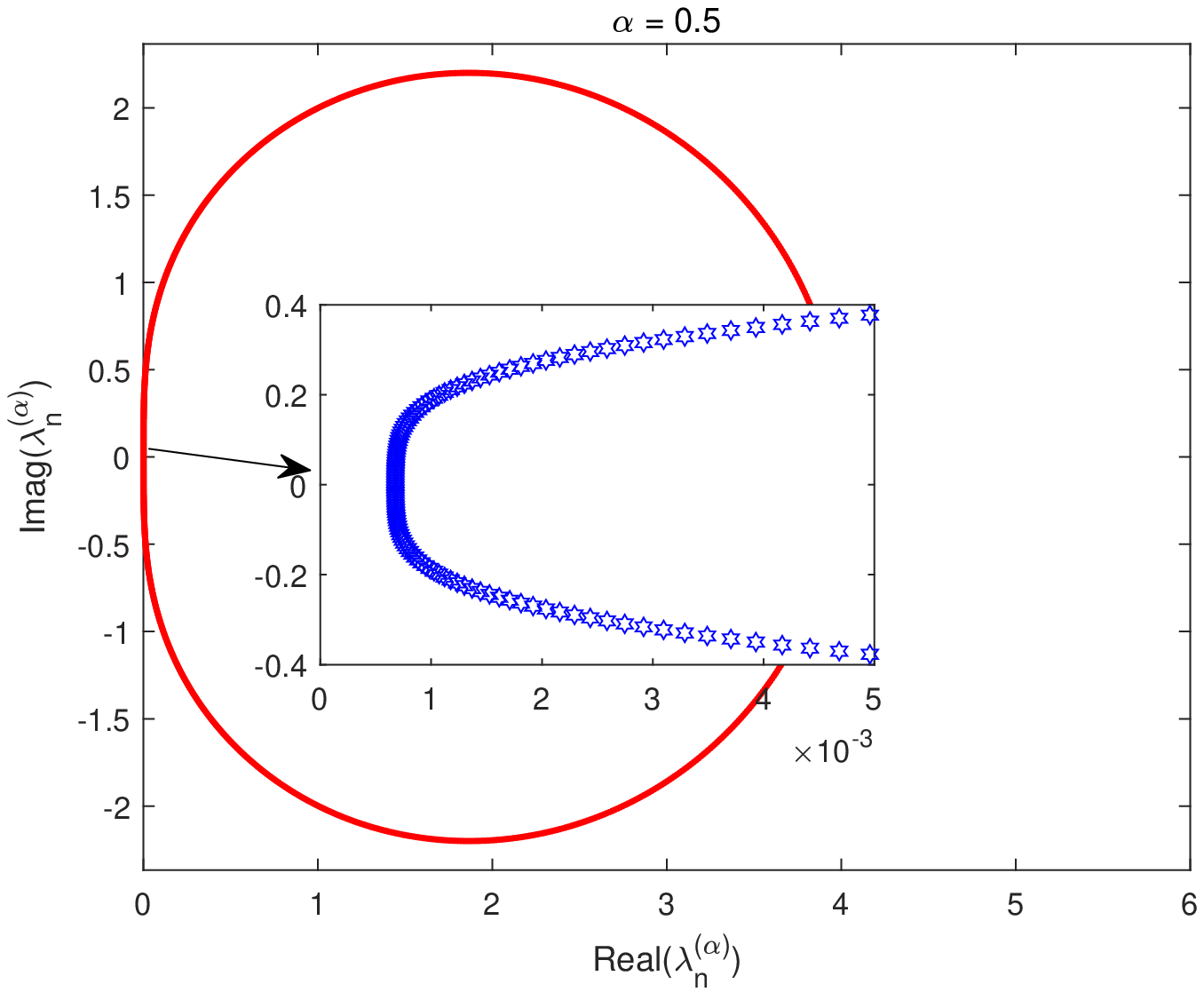}
\includegraphics[width=2.137in,height=1.9in]{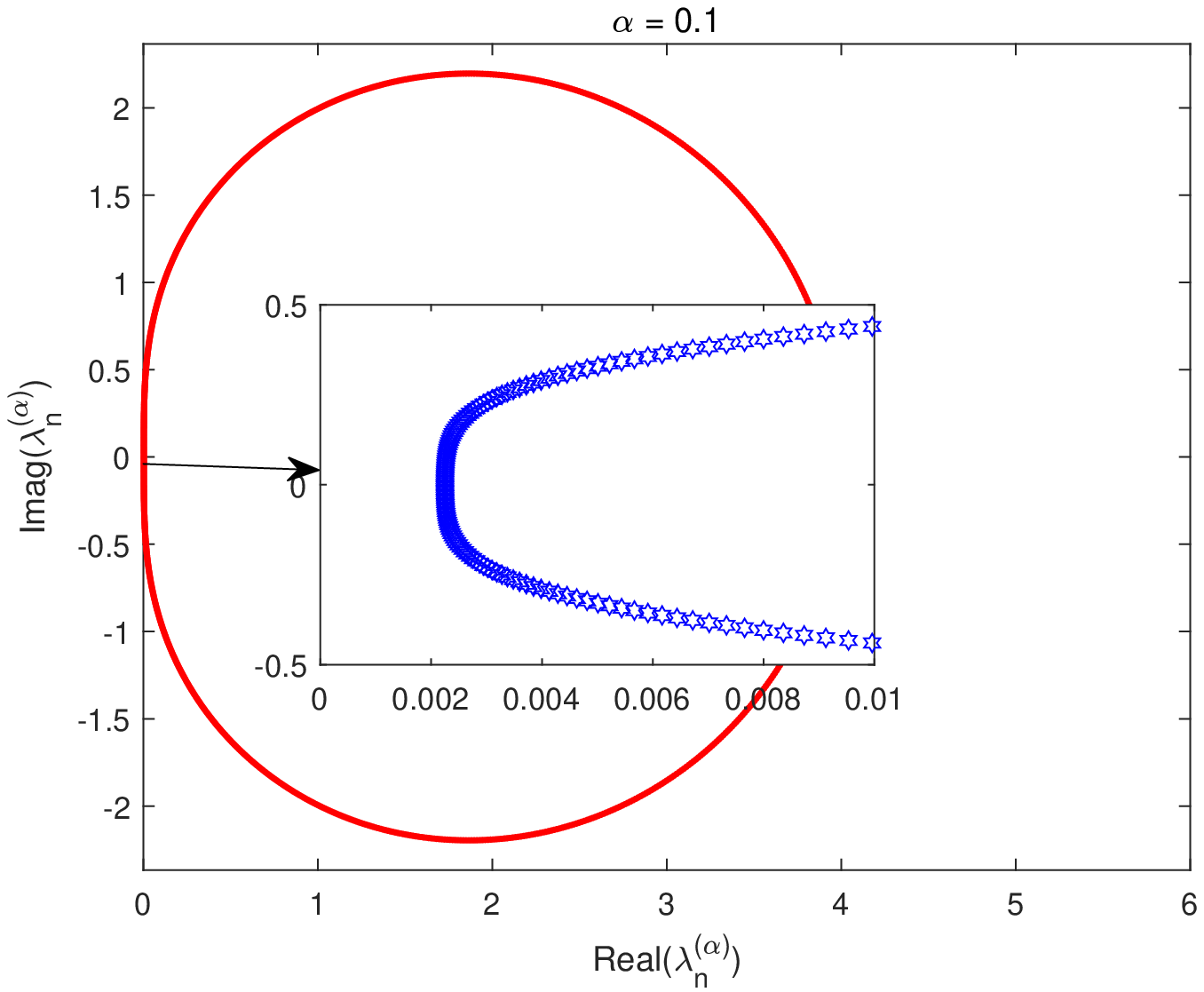}
\caption{The plot of complex quantities $\{\lambda^{(\alpha)}_n\}^{
N_t}_{n=1}$ with different $\alpha$'s.}
\label{fig3.1}
\end{figure}

In practice, we often choose $\alpha\in(0,1)$ \cite{Lin20}; then it holds that $
\mathcal{R}e(\lambda^{(\alpha)}_n) > 0$ according to Proposition \ref{pro3.1}. In Fig.~\ref{fig3.1}, we plot the complex
quantities $\{\lambda^{(\alpha)}_n\}^{N_t}_{n=1}$ on the complex plane. We see that for $\alpha\in(0,1)$, it holds $\mathcal{R}e(\lambda^{(\alpha)}_n)> 0$ and consequently
all the linear systems involved at Step-(b) are positive definite, thus not difficult to solve~\cite{Banjai12}. As shown in~\cite{Wu2018}, a multigrid method using Richardson iterations as a smoother with an optimized choice of the damping parameter is very effective. By the way, a detailed convergence analysis with multigrid precondioning in the context of the $\tau$ algebra considered in this paper is presented in~\cite{Domase}.

In addition, it is not hard to derive the following result.
\begin{theorem}
Let $\alpha\in(0,1]$ and $\mathcal{P}^{-1}_{\alpha}\mathcal{A} = \mathcal{I}
+ \mathcal{L}$, where $\mathcal{I}$ is an identity matrix of order $N_t
N_s$, then it follows that ${\rm rank}(\mathcal{L}) = 2N_t$.
\label{thm3.1}
\end{theorem}
\noindent\textbf{Proof}. We consider that
\begin{equation}
\begin{split}
\mathcal{P}^{-1}_{\alpha}\mathcal{A} & = \mathcal{P}^{-1}_{\alpha}[\mathcal{P}_{\alpha} +
(\mathcal{A} - \mathcal{P}_{\alpha})]\\
& = \mathcal{I} + \mathcal{P}^{-1}_{\alpha}[(C - C_{\alpha})\otimes I_s]\\
&\triangleq \mathcal{I} + \mathcal{L},
\end{split}
\end{equation}
where ${\rm rank}(\mathcal{L}) = {\rm rank}\big(\mathcal{P}^{-1}_{\alpha}[(C -
C_{\alpha})\otimes I_s]\big) = {\rm rank}\big((C - C_{\alpha})\otimes I_s
\big) = 2N_t$, which completes the proof. \hfill$\Box$

Theorem \ref{thm3.1} implies that the generalized minimal residual
(GMRES) method by Saad and Schultz~\cite{Saad86} can compute the exact solution
of the preconditioned system $\mathcal{P}^{-1}_{\alpha}\mathcal{A}$ (if it is
diagonalizable) in at most $2N_t + 1$ iterations \cite{Saad86,Gu2015,McDona}.
Although it should be noted that such a convergence behavior of the
preconditioned GMRES is not sharp
when $N_t$ is not small, the result can give useful clues on the effectiveness
of the preconditioner $\mathcal{P}_{\alpha}$ to approximate the all-at-once matrix
$\mathcal{A}$. Moreover, in Section \ref{sec4} we will show numerically that most of the eigenvalues of the preconditioned matrix $\mathcal{P}^{-1}_{\alpha}\mathcal{A}$ actually cluster at point 1 of the spectrum, supporting the theoretical findings of this section on the good potential of the proposed preconditioner $\mathcal{P}_{\alpha}$ to accelerate the iterative solution of $\mathcal{P}^{-1}_{\alpha}\mathcal{A}$.


\begin{remark}
If we set $\alpha = 1$, the  preconditioner $\mathcal{P}_{\alpha}$ reduces
to the BC preconditioner $\mathcal{P}_1 = C_1\otimes I_s - \tau I_t\otimes A$,
where $C_1$ is a matrix of $C_{\alpha}$ in (\ref{eq1.11}) with $\alpha = 1$.
Moreover, the matrix decomposition (\ref{eq1.7x}), Proposition \ref{pro3.1},
and Theorem \ref{thm3.1} are available for the BC preconditioner $\mathcal{P}_1$,
refer e.g., to \cite{McDona,Lin20} for details;
\end{remark}
\begin{remark}
We remark that the invertibility of $\mathcal{P}_{\alpha}$ (or $\mathcal{P}_1$)
depends completely on the matrix $\lambda^{(\alpha)}_{n} I_s - \tau A$. However,
if we choose $\alpha \rightarrow 1$, the real parts of $\{\lambda^{(\alpha)}_n\}^{N_t}_{n=1}$
will be very close to zero, cf. Proposition \ref{pro3.1}. Meanwhile, if the diffusion coefficient $\kappa_{\gamma}$ is very small, the eigenvalues of $A$ will be increasingly close to zero. This fact will make the matrices $\lambda^{(\alpha)}_{n} I_s - \tau A$ potentially very ill-conditioned (even singular\footnote{If the model problem (\ref{eq1.1}) incorporates Neumann boundary condition and we set $\alpha = 1$, the matrices $\lambda^{(\alpha)}_{n} I_s - \tau A$ will be singular because of the singular matrix $A$.}). Under this circumstance, the BC preconditioner $\mathcal{P}_1$ should be a bad preconditioner (see \cite{Lin20} and Section \ref{sec4}) while the preconditioner $\mathcal{P}_{\alpha}$ with $\alpha\in(0,1)$ should be preferred in a practical implementation;
\label{rem3.2}
\end{remark}
\begin{remark}
If one wants to improve the proposed preconditioner further, then the
polynomial preconditioner that has potential to reduce communication costs in Krylov
subspace solvers \cite{Saad86} can be derived:
\begin{equation}
\begin{split}
\mathcal{A}^{-1} &= \left[\mathcal{I} - \mathcal{P}^{-1}_{\alpha}\mathcal{S}_{\alpha}
\right]^{-1}\mathcal{P}^{-1}_{\alpha}\\
& = \sum^{\infty}_{k=1}(\mathcal{P}^{-1}_{\alpha}\mathcal{S}_{\alpha})^k\mathcal{P}^{-1}_{\alpha},
\end{split}
\end{equation}
where $\mathcal{S}_{\alpha} = (C_{\alpha}- C)\otimes I_x$ and assume that $\|\mathcal{P}^{-1}_{\alpha}\mathcal{S}_{\alpha}\| < 1$, then the $m$-step polynomial
preconditioner $\mathcal{P}_{\alpha}(m) = \sum^{m}_{k=0}(\mathcal{P}^{-1}_{\alpha}\mathcal{
S}_{\alpha})^k\mathcal{P}^{-1}_{\alpha}$ will be available and inherently parallel. Moreover, such a preconditioner
can be easily adapted for the all-at-once systems arising from the given non-uniform
temporal discretizations; refer to \cite{Goddard} for a short discussion.
\end{remark}
%
%
%
%
\subsection{Efficient implementation of the PinT preconditioner}
In this subsection, we discuss on how to implement the Krylov subspace method
more efficiently for solving the left (or right) preconditioned system $\mathcal{P
}^{-1}_{\alpha}\mathcal{A}{\bm U} = \mathcal{P}^{-1}_{\alpha}{\bm F}$ (or $\mathcal{A}
\mathcal{P}^{-1}_{\alpha}(\mathcal{P}_{\alpha}{\bm U}) = {\bm F}$), $\alpha\in(0,1]$. We recall that the kernel operations of a Krylov subspace algorithm are the matrix-vector products $
\mathcal{P}^{-1}_{\alpha}{\bm v}$ and $\mathcal{A}{\bm v}$ for some given vector
${\bm v}\in\mathbb{R}^{N_tN_s\times 1}$. First, we present a fast implementation for
computing the second matrix-vector product.
\begin{equation}
\begin{split}
\mathcal{A}{\bm v} & = (C\otimes I_s){\bm v} - (\tau I_t\otimes A){\bm v}\\
&={\rm vec}(VC^{T}) - \tau\cdot{\rm vec}(AV),
\end{split}
\end{equation}
where the operation ${\bm v} = {\rm vec}(V)$ and ${\rm vec}(X)$ denote the vectorization of the matrix $X$ obtained by stacking the columns of $X$ into a single column vector. The
first term can be calculated directly in $\mathcal{O}(N_tN_s)$ storage and operations
due to the sparse matrix $C^T$. Owing to the Toeplitz structure of matrix $A$, the second term can
be evaluated in $\mathcal{O}(N_tN_s\log N_s)$ operations and $\mathcal{O}(N_tN_s)$
storage. This means that the computation of $\mathcal{A}{\bm v}$ requires the same
operations and storage of $(I_t\otimes A){\bm v}$. Clearly, the algorithm complexity of $\mathcal{A}{\bm v}$ can be further alleviated by a parallel implementation.

Next, we focus on the fast computation of $\mathcal{P}^{-1}_{\alpha}{\bm v}$ for a given vector ${\bm v}$. According to Step-(a)--Step-(c) described in Section \ref{sec3.1}, the first and third steps of the implementation of $\mathcal{P}_{\alpha}$ require $\mathcal{O}(N_tN_s\log N_s)$ operations and $\mathcal{O}(N_tN_s)$ storage.
If the matrix $A$ can be diagonalized (or approximated) via fast DSTs, Step-(b) can be directly solved with $\mathcal{O}(N_tN_s\log N_s)$ operations and $\mathcal{O}(N_tN_s)$ storage, except for the different complex shifts.

Similar to \cite[pp. 10-11]{Lin20}, it is worth noting that {\em only half of $N_t$
shifted linear systems in Step-(b) need to be solved}. From Step-(a), we know
that the right hand sides in Step-(b) can be expressed as
\begin{equation}
{\bm z}_{1,n} = \frac{1}{\sqrt{N_t}}\sum^{2}_{j = 0}
\alpha^{\frac{j}{N_t}}\theta^{(n-1)j}{\bm v}_{n},~~ {\bm v} = [{\bm v}^{\top}_{1},
{\bm v}^{\top}_{2},\cdots,{\bm v}^{\top}_{N_t}]^{\top},~~ n = 1,2,\cdots,N_t.
\label{eq3.10x}
\end{equation}
We recall the expression of coefficient matrices in Step-(b):
\begin{equation}
B_{n-1} := \lambda^{(\alpha)}_{n} I_s - \tau A = \left(\sum^{2}_{j = 0}r_j
\alpha^{\frac{j}{N_t}}\theta^{(n-1)j}\right)I_s - \tau A,\quad n = 1,2,\cdots,N_t.
\end{equation}

Let ``${\rm conj}(\cdot)$" denote conjugate of a matrix or a vector. Then,
\begin{equation}
\begin{cases}
{\rm conj}({\bm z}_{1,n}) = \frac{1}{\sqrt{N_t}}\sum\limits^{2}_{j = 0}
\alpha^{\frac{j}{N_t}}\theta^{[N_t-(n-1)]j}{\bm v}_{n} = {\bm z}_{1,N_t - n + 2},\\
{\rm conj}(B_{n-1}) = \left(\sum\limits^{2}_{j = 0}r_j
\alpha^{\frac{j}{N_t}}\theta^{-(n-1)j}\right)I_s - \tau A = \left(\sum\limits^{2}_{j = 0}r_j
\alpha^{\frac{j}{N_t}}\theta^{[N_t-(n-1)]j}\right)I_s - \tau A = B_{N_t - {\cre{n}} + 1},& n = 2,3,\cdots,N_t.
\end{cases}
\label{eq3.12x}
\end{equation}
That means the unknowns in Step-(b) hold equalities: ${\bm z}_{2,n} = {\rm conj}({\bm z}_{2,N_t-n+2})$
for $n = 2, 3,\cdots, N_t$. Therefore, only the first $\lfloor(N_t+1)/2\rfloor$ multi-shifted linear systems in Step-(b) need to be solved, thus the number of core processors required in the practical parallel implementation is significantly reduced.

Finally, we discuss the efficient solution of the sequence of shifted linear systems in Step-(b). Since $A$ is a real negative definite Toeplitz matrix, it can be approximated efficiently by a $\tau$-matrix $\mathscr{T}(A)$~\cite{Bini90} that is a real symmetric matrix and can be diagonalized as follows:
\begin{equation}
\mathscr{T}(A) = Q^{\top}_{N_s}\Lambda_{N_s} Q_{N_s},\quad Q_{N_s} = \sqrt{\frac{2}{N_s+1}}
\sin\left(\frac{\pi ij}{N_s + 1}\right),~1\leq i,j\leq N_s,\quad\Lambda_{N_s} = {\rm diag}(\sigma_1,
\sigma,\cdots,\sigma_{N_s}).
\label{eq3.9}
\end{equation}
From the relations $Q^{\top}_{N_s}Q_{N_s} = I_s$ and $Q_{N_s}\mathscr{T}(A) = \Lambda_{N_s}Q_{N_s}$, we have
\begin{equation}
\sigma_i = \frac{\sum^{N_s}_{\ell = 1}a_{\ell}\sin\left(\frac{\pi ij}
{N_s + 1}\right)}{\sin\left(\frac{\pi i}{N_s + 1}\right)},
\end{equation}
where the vector $[a_1,~a_2,\cdots,a_{N_s}]^{\top}$ is the first column of $A$. The set of shifted linear systems at Step-(b) is reformulated as the following sequence:
\begin{equation}
\big[\lambda^{(\alpha)}_{n} I_s - \tau\mathscr{T}(A)\big]{\bm z}_{2,n} =
Q^{\top}_{N_s}\big[\lambda^{(\alpha)}_{n} I_s - \tau\Lambda_{N_s}\big]
Q_{N_s} = {\bm z}_{1,n},\quad  ~n=1,2,\dots, N_t,
\label{eq3.11}
\end{equation}
that can be solved by fast DSTs without storing any dense matrix. In addition, since the DST only involves real operations, so the above
preconditioner matrix diagonalization does not affect Eqs. (\ref{eq3.10x})--(\ref{eq3.12x}) at all. Overall, the fast computation of $P^{-1}_{\alpha}{\bm v}$ using Eq. (\ref{eq3.11}) requires $\mathcal{O}(N_tN_s\log N_s)$ operations
and $\mathcal{O}(N_tN_s)$ memory units.
\begin{remark}
It is worthwhile noting that such an implementation is also suitable for
solving the all-at-once system that arises from the spatial discretizations
using the compact finite difference, finite element and spectral methods by
only substituting in Eq.~(\ref{eq3.11}) the identity matrix $I_s$ with the
mass matrix. Fortunately, for one-dimensional problems, the mass matrix is a SPD Toeplitz tridiagonal matrix~\cite{Lei17,Ran2019,Xu2019}, which can be
naturally diagonalized by fast DSTs. For high-dimensional models, the mass
matrix will be a SPD Toeplitz block tridiagonal matrix in the Kronecker product form, thus it can be still diagonalized via fast DSTs; refer, e.g., to \cite{Hu2019};
\end{remark}
\begin{remark}
Finally, it is essential to mention that although the algorithm and theoretical
analyses are presented for one-dimensional model (\ref{eq1.1}), they can be naturally
extended for two-dimensional model
(refer to the next section for a discussion) because the above algorithm and theoretical
analyses are mainly based on the temporal discretization.

%
\end{remark}

\section{Numerical experiments}
\label{sec4}
In this section, we present the performance of the proposed $\alpha$-circulant-based
preconditioners for solving  some examples of one- and two-dimensional model problems
(\ref{eq1.1}). In all our numerical experiments
reported in this section, following the guidelines given in \cite{Lin20} we set $\alpha=\min\{
0.5, 0.5\tau\}$ for the preconditioners $\mathcal{P}_{\alpha}$ and $\mathcal{P}_1$. All our  experiments are performed in MATLAB R2017a on a Gigabyte workstation equipped with Intel(R) Core(TM)
i9-7900X CPU @3.3GHz, 32GB RAM running the Windows 10 Education operating system (64bit version) using double precision floating point arithmetic (with machine epsilon equal to $10^{-16}$). The adaptive
Simpler GMRES (shortly referred to as Ad-SGMRES) method~\cite{Pavel10} is employed to solve the right-preconditioned systems in Example~1, while the BiCGSTAB method~\cite[pp. 244-247]{Saad86} method is applied to the left-preconditioned systems in Example 2 using the built-in function available in MATLAB\footnote{In this case, since Ad-SGMRES
still requires large amounts of storage due to the orthogonalization process, we have also
used the BiCGSTAB method as an alternative iterative method for solving non-symmetric
systems.}. The tolerance for the stopping criterion in both algorithms is set as $\|{\bm r}_k\|_2/\|
{\bm r}_0\|_2 <{\rm tol}=10^{-9}$, where ${\bm r}_k$ is the residual vector at $k$th Ad-SGMRES or BiCGSTAB iteration. The iterations are always started from the zero vector. All timings shown in the tables are obtained by averaging over 20 runs\footnote{The code in MATLAB will be available at \url{https://github.com/Hsien-Ming-Ku/Group-of-FDEs}.}.

In the tables, the quantity `\texttt{Iter}' represents the iteration number of Ad-SGMRES or BiCGSTAB, `DoF' is the number of degrees of freedom, or unknowns, and `\texttt{CPU}' is the computational time expressed in seconds. The 2-norm of the true relative residual (called in the tables  \texttt{TRR}) is defined as
$\mathtt{TRR} = ||{\bm F} - \mathcal{A}{\bm U}_k\|_2/\|{\bm F}\|_2$, and the numerical error (\texttt{Err}) between the approximate and the exact solutions at the final time level reads $\|{\bm u}^{*} - {\bm u}^{N_t}\|_{\infty}$, where ${\bm U}_k$ is the approximate solution when the preconditioned iterative solvers terminate and ${\bm u}^{*}$ is the exact solution on the mesh. These notations are adopted throughout this section.
\vspace{1.5mm}

\noindent\textbf{Example 1.} The first example is a Riesz fractional diffusion
equation (\ref{eq1.1}) with coefficients $(a,b) = (0,1)$, $T = 1$, $\kappa_{\gamma}=0.01$ and $\phi(x) = 15(1 + \gamma/4)x^3(1 - x)^3$. The source term is given by
\begin{equation*}
\begin{split}
f(x,t) & = 15(1 + \gamma/4)e^tx^3(1 - x)^3 + \frac{15(1 + \gamma/4)\kappa_{\gamma}e^t}{2\cos(\gamma\pi/2)}
\Bigg[\frac{\Gamma(4)}{\Gamma(4 - \gamma)}\big(x^{3 - \gamma} + (1 - x)^{3
- \gamma}\big) - \frac{3\Gamma(5)}{\Gamma(5 - \gamma)}\big(x^{4 - \gamma}
~+ \\
&\quad~(1 - x)^{4 - \gamma}\big) +\frac{3\Gamma(6)}{\Gamma(6 - \gamma)}\big(x^{5 - \gamma} + (1 - x)^{5 -
\gamma}\big) - \frac{\Gamma(7)} {\Gamma(7 - \gamma)}\big(x^{6 - \gamma} +
(1 - x)^{6 - \gamma}\big)
\Bigg].
\end{split}
\end{equation*}
The exact solution is known and it reads as $u(x,t) = 15(1 + \gamma/4)e^tx^3(1 - x)^3$. The results of our numerical experiments with the Ad-SGMRES method preconditioned by $\mathcal{P}_{\alpha}$ and $\mathcal{P}_1$ for solving the all-at-once
discretized systems (\ref{eq12}) are reported in Tables \ref{table:nonlin}--\ref{tab3}.

\begin{figure}[!htpb]
\centering
\includegraphics[width=3.12in,height=2.65in]{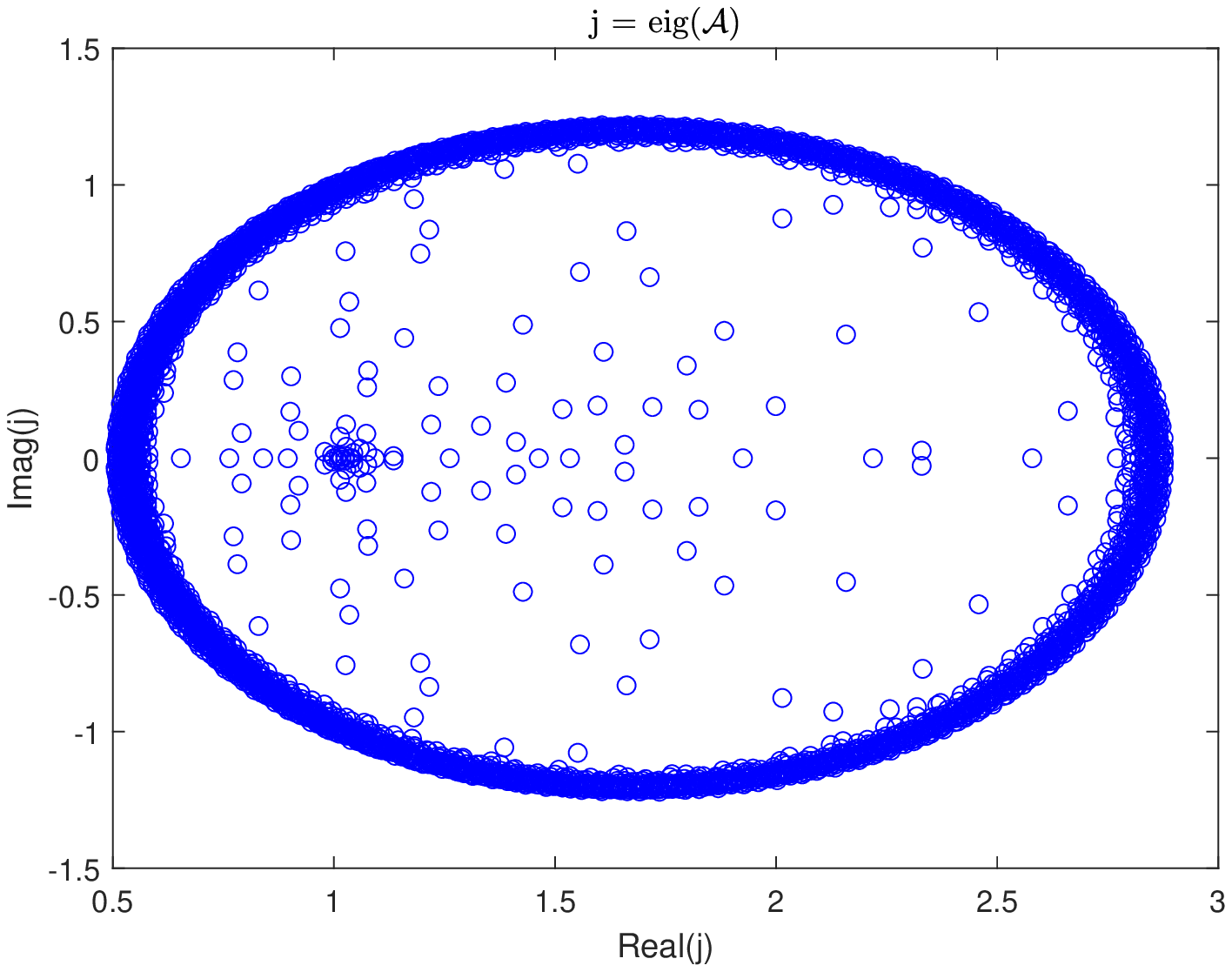} 
\includegraphics[width=3.12in,height=2.65in]{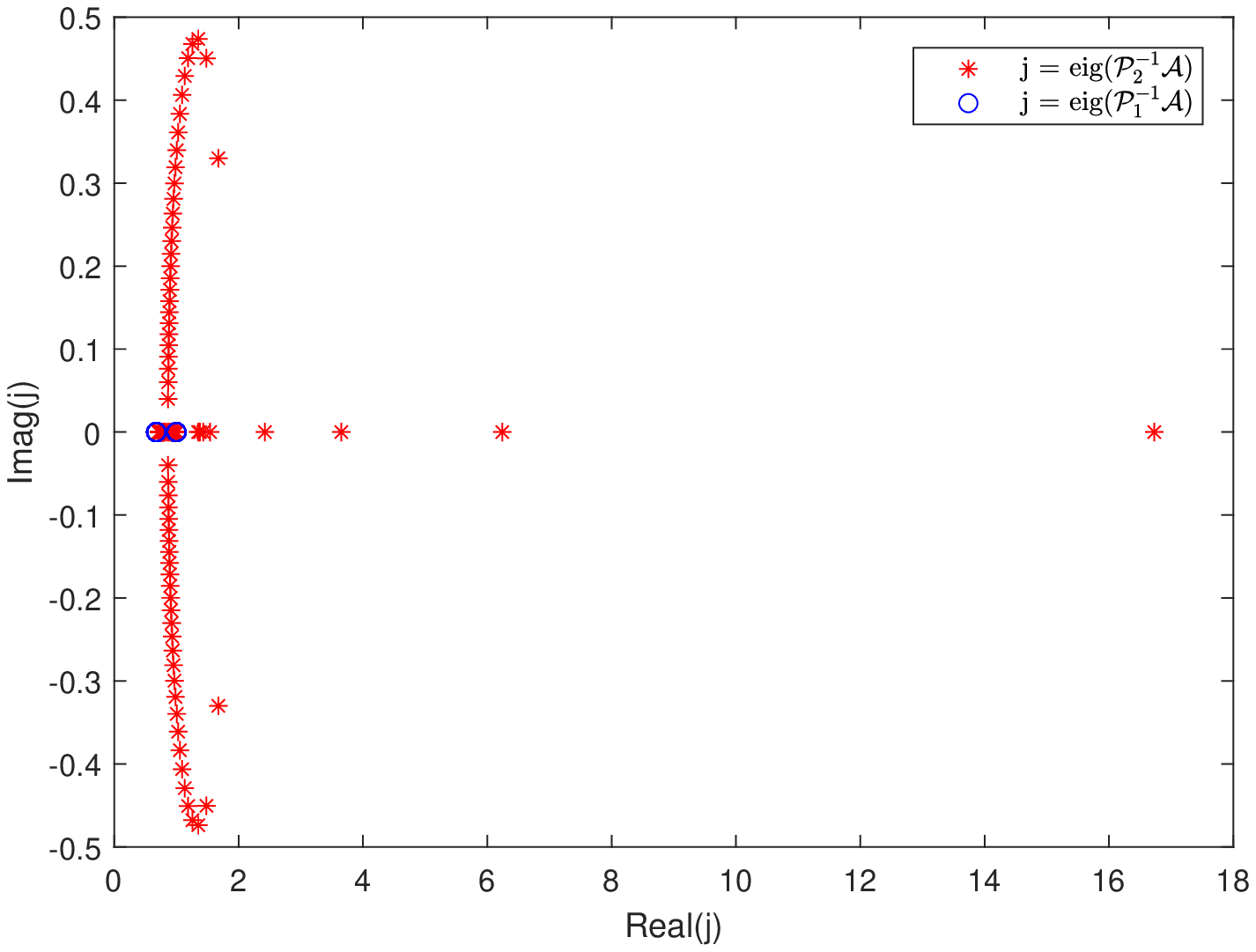}
\includegraphics[width=3.12in,height=2.65in]{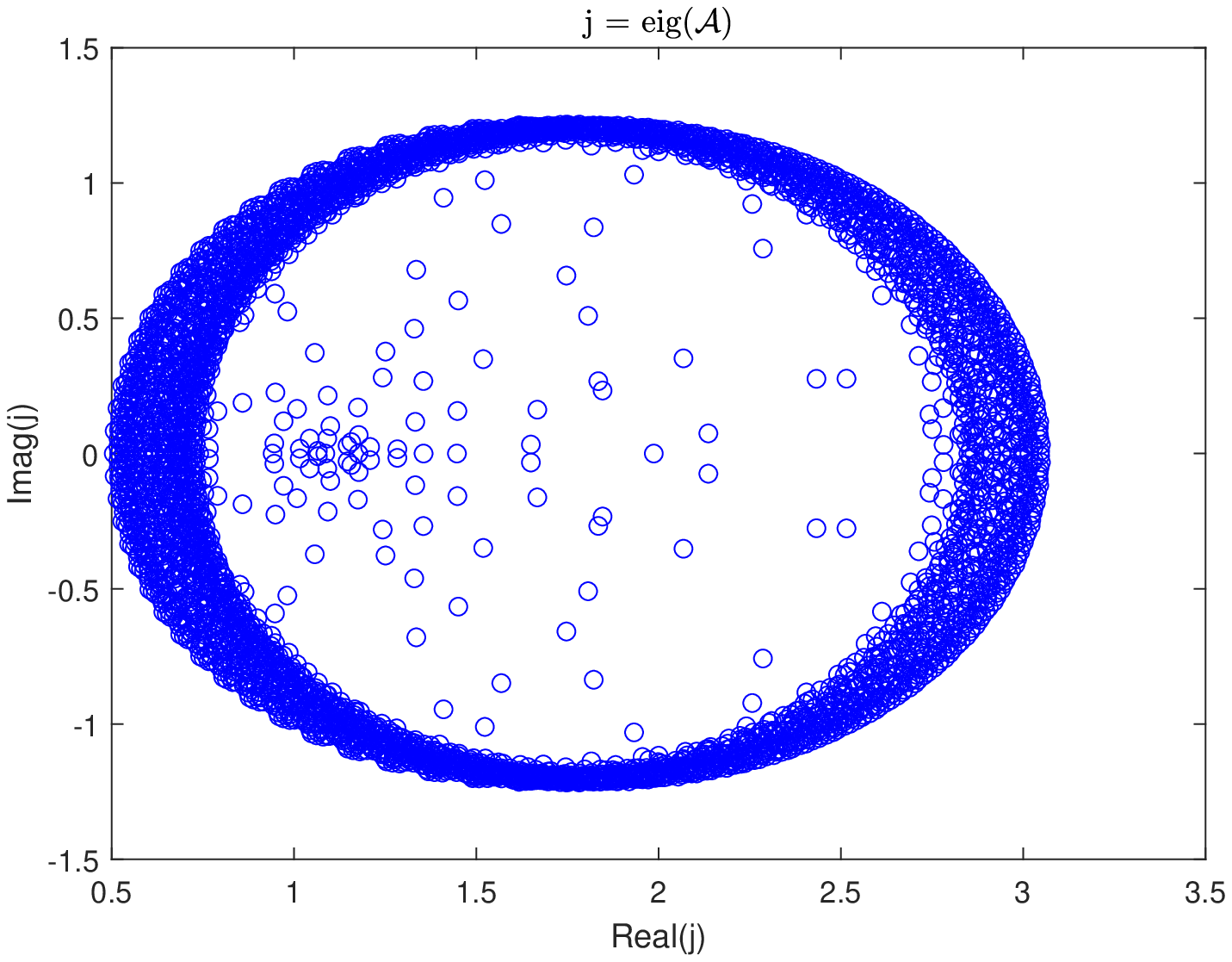}
\includegraphics[width=3.12in,height=2.65in]{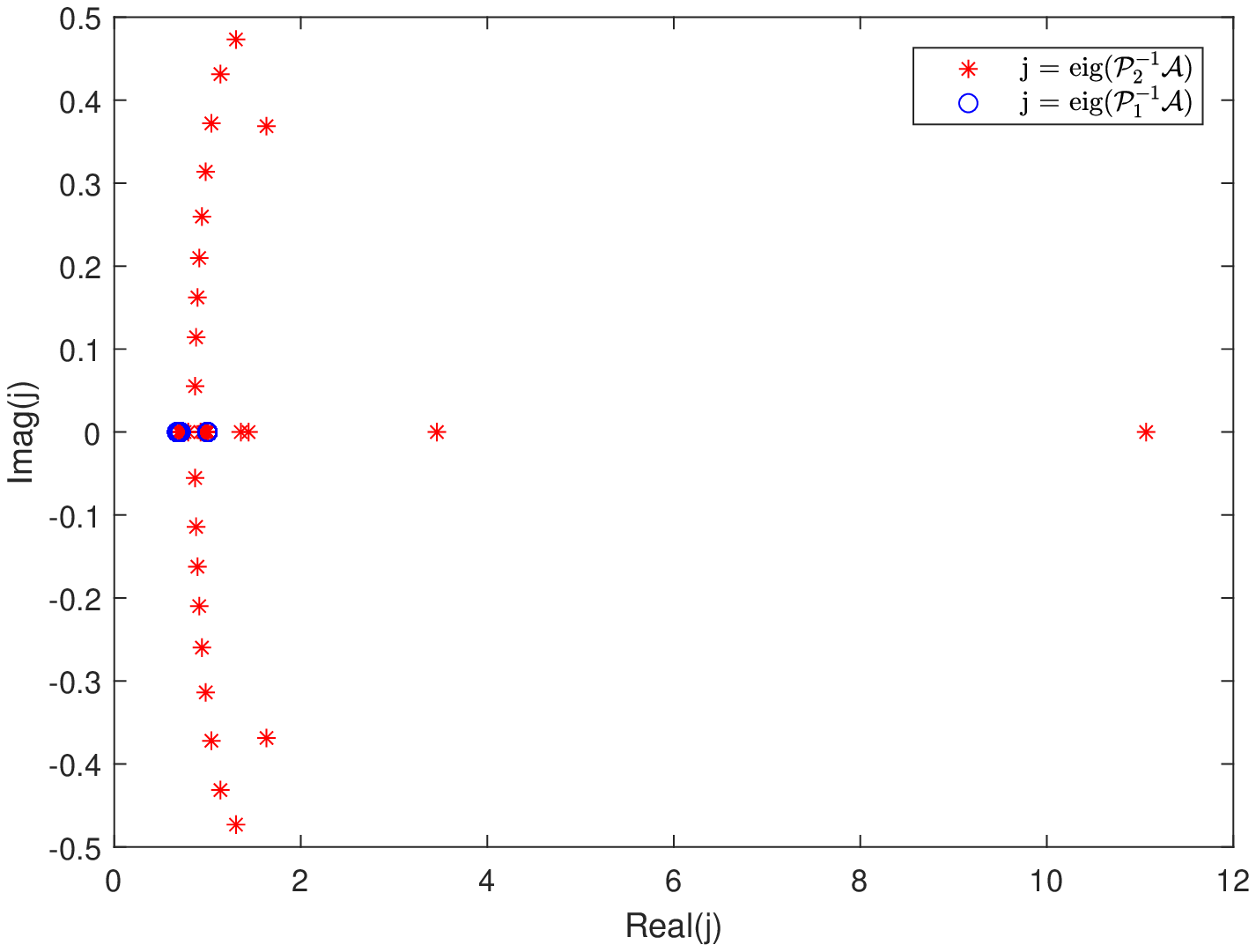}
\includegraphics[width=3.12in,height=2.65in]{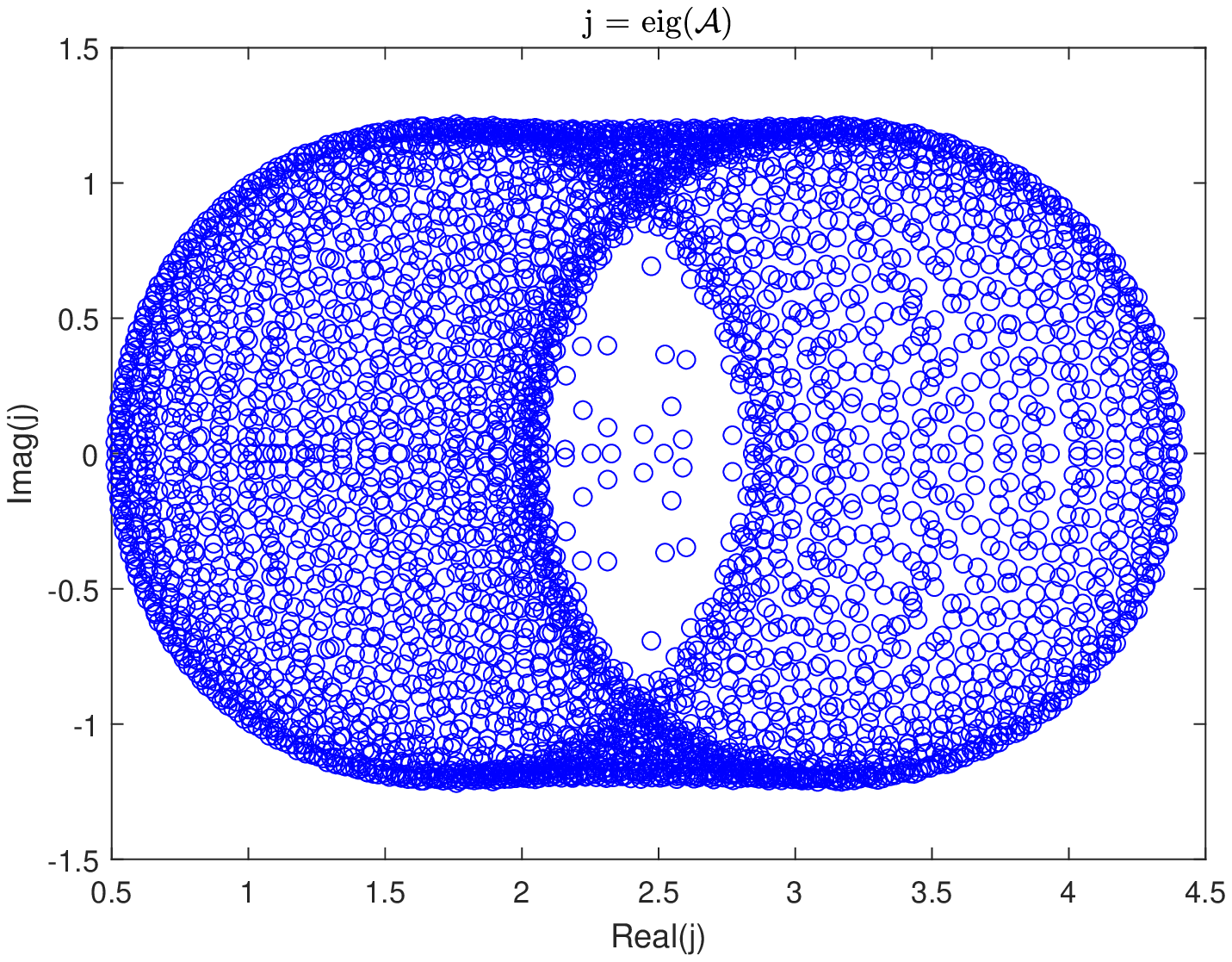}
\includegraphics[width=3.12in,height=2.65in]{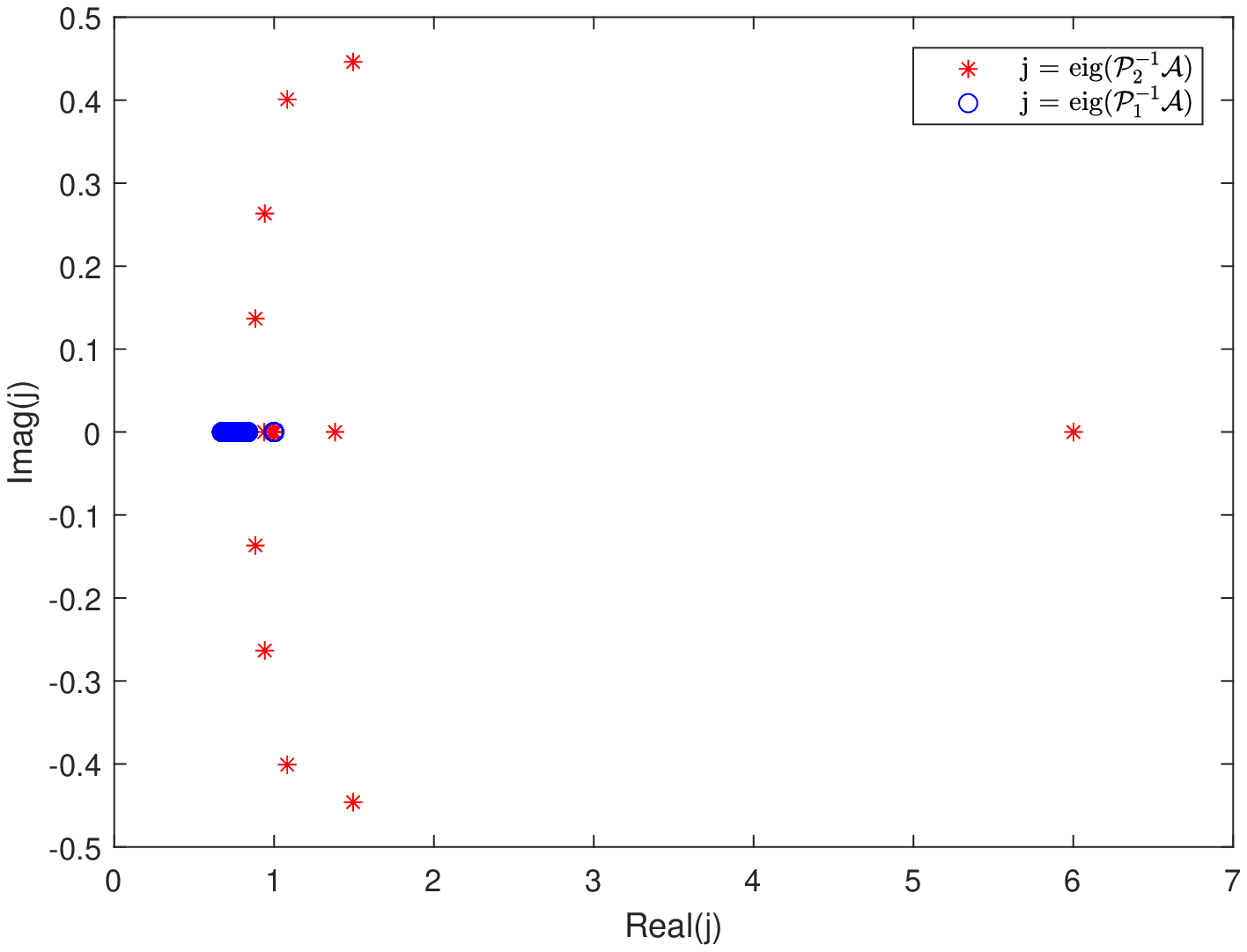}
\caption{The eigenvalue distributions of the matrix $\mathcal{A}$
and preconditioned matrices $\mathcal{P}^{-1}_{\alpha}\mathcal{A},~
\mathcal{P}^{-1}_1\mathcal{A}$ with
$\gamma = 1.2,~1.5,~1.9$ and $h = \tau = 1/64$ in Example 1.}
\label{fig4.1}
\end{figure}

\begin{table}[!htpb]
\begin{center}
\caption{Numerical results of GMRES with two different preconditioners on Example 1 with
$\gamma = 1.2$.}
\label{table:nonlin}
\vspace{1mm}
\begin{tabular}{crrccrccccc}
\hline &&&\multicolumn{4}{c}{$\mathcal{P}_{\alpha}$} &\multicolumn{4}{c}{$\mathcal{P}_1$}\\
[-2pt]\cmidrule(l{0.7em}r{0.7em}){4-7} \cmidrule(l{0.7em}r{0.6em}){8-11}\\[-11pt]
$N_t$ &$h$ &DoF &\texttt{Iter} &\texttt{CPU} &\texttt{TRR} &\texttt{Err}
&\texttt{Iter} &\texttt{CPU} &\texttt{TRR} &\texttt{Err} \\
\hline
$2^6$    &1/128    &8,128     &7 &0.033 &-9.794  &9.7599e-5 &19 &0.085 &-9.412 &9.7599e-5 \\
         &1/256    &16,320    &7 &0.054 &-9.257  &9.4838e-5 &19 &0.117 &-9.347 &9.4838e-5  \\
         &1/512    &32,704    &8 &0.116 &-10.017 &9.4147e-5 &19 &0.247 &-9.535 &9.4147e-5  \\
         &1/1024   &65,472    &8 &0.156 &-9.537  &9.3974e-5 &19 &0.359 &-9.163 &9.3975e-5  \\
\hline
$2^8$    &1/128    &32,512    &7 &0.107 &-10.366 &9.5721e-6 &19 &0.258 &-9.378 &9.5722e-6  \\
         &1/256    &65,280    &7 &0.158 &-9.580  &6.8110e-6 &19 &0.385 &-9.318 &6.8111e-6  \\
         &1/512    &130,816   &7 &0.445 &-9.020  &6.1205e-6 &19 &1.236 &-9.248 &6.1208e-6	 \\
         &1/1024   &261,888   &8 &0.787 &-9.763  &5.9481e-6 &19 &1.958 &-9.156 &5.9482e-6  \\
\hline
$2^{10}$ &1/128    &130,048   &6 &0.323 &-9.020  &5.0121e-6 &19 &0.964 &-9.368 &5.0138e-6	 \\
         &1/256    &261,120   &7 &0.756 &-9.895  &1.2888e-6 &19 &2.123 &-9.309 &1.2890e-6 \\
         &1/512    &523,264   &7 &1.715 &-9.326  &5.9821e-7 &19 &4.749 &-9.242 &5.9870e-7 \\
         &1/1024   &1,047,552 &8 &3.245 &-10.056 &4.2607e-7 &19 &8.103 &-9.153 &4.2613e-7 \\
\hline
\end{tabular}
\end{center}
\end{table}
\begin{table}[t]
\begin{center}
\caption{Numerical results of GMRES with two different preconditioners on Example 1
with $\gamma = 1.5$.}
\vspace{1mm}
\begin{tabular}{crrccccccrc}
\hline &&&\multicolumn{4}{c}{$\mathcal{P}_{\alpha}$} &\multicolumn{4}{c}{$\mathcal{P}_1$}\\
[-2pt]\cmidrule(l{0.7em}r{0.7em}){4-7} \cmidrule(l{0.7em}r{0.6em}){8-11}\\[-11pt]
$N_t$ &$h$ &DoF &\texttt{Iter} &\texttt{CPU} &\texttt{TRR} &\texttt{Err}
&\texttt{Iter} &\texttt{CPU} &\texttt{TRR} &\texttt{Err} \\
\hline
$2^6$    &1/128    &8,128     &8 &0.038 &-9.892  &1.0514e-4 &15 &0.085 &-9.044 &1.0515e-4 \\
         &1/256    &16,320    &8 &0.054 &-9.755  &9.8789e-5 &15 &0.095 &-9.176 &9.8789e-5 \\
         &1/512    &32,704    &8 &0.112 &-9.697  &9.7199e-5 &15 &0.198 &-9.513 &9.7199e-5 \\
         &1/1024   &65,472    &8 &0.157 &-9.349  &9.6802e-5 &16 &0.308 &-9.563 &9.6802e-5 \\
\hline
$2^8$    &1/128    &32,512    &7 &0.112 &-9.545  &1.4536e-5 &16 &0.225 &-9.831 &1.4536e-5 \\
         &1/256    &65,280    &7 &0.158 &-9.103  &8.1809e-6 &15 &0.319 &-9.323 &8.1810e-6	\\
         &1/512    &130,816   &8 &0.515 &-9.925  &6.5922e-6 &15 &0.979 &-9.265 &6.5922e-6 \\
         &1/1024   &261,888   &8 &0.829 &-9.587  &6.1950e-6 &16 &1.649 &-9.605 &6.1950e-6 \\
\hline
$2^{10}$ &1/128    &130,048   &7 &0.392 &-9.982  &1.3161e-5 &15 &0.827 &-9.008 &1.3162e-5 \\
         &1/256    &261,120   &7 &0.779 &-9.423  &3.2696e-6 &15 &1.637 &-9.184 &3.2692e-6	\\
         &1/512    &523,264   &7 &1.744 &-9.049  &9.0813e-7 &15 &3.769 &-9.381 &9.0882e-7 \\
         &1/1024   &1,047,552 &8 &3.221 &-9.878  &5.1171e-7 &16 &6.658 &-9.631 &5.1160e-7 \\
\hline
\end{tabular}
\end{center}
\label{tab2}
\end{table}
\begin{table}[!htpb]
\begin{center}
\caption{Numerical results of GMRES with two different preconditioners on Example 1 with
$\gamma = 1.9$.}
\label{tab3}
\vspace{1mm}
\begin{tabular}{crrcccccccc}
\hline &&&\multicolumn{4}{c}{$\mathcal{P}_{\alpha}$} &\multicolumn{4}{c}{$\mathcal{P}_1$}\\
[-2pt]\cmidrule(l{0.7em}r{0.7em}){4-7} \cmidrule(l{0.7em}r{0.6em}){8-11}\\[-11pt]
$N_t$ &$h$ &DoF &\texttt{Iter} &\texttt{CPU} &\texttt{TRR} &\texttt{Err}
&\texttt{Iter} &\texttt{CPU} &\texttt{TRR} &\texttt{Err} \\
\hline
$2^6$    &1/128    &8,128     &7  &0.034 &-9.255  &1.2052e-4 &11 &0.052 &-9.305 &1.2052e-4 \\
         &1/256    &16,320    &7  &0.047 &-9.101  &1.0303e-4 &11 &0.070 &-9.195 &1.0303e-4 \\
         &1/512    &32,704    &8  &0.113 &-10.012 &9.8653e-5 &11 &0.151 &-9.123 &9.8653e-5 \\
         &1/1024   &65,472    &8  &0.159 &-10.024 &9.7559e-5 &11 &0.225 &-9.058 &9.7559e-5 \\
\hline
$2^8$    &1/128    &32,512    &7  &0.112 &-9.831  &3.8671e-5 &11 &0.153 &-9.408 &3.8671e-5 \\
         &1/256    &65,280    &7  &0.155 &-9.500  &1.1924e-5 &11 &0.235 &-9.313 &1.1924e-5 \\
         &1/512    &130,816   &7  &0.456 &-9.287  &7.5514e-6 &11 &0.690 &-9.255 &7.5518e-6 \\
         &1/1024   &261,888   &7  &0.696 &-9.261  &6.4585e-6 &11 &1.093 &-9.212 &6.4587e-6 \\
\hline
$2^{10}$ &1/128    &130,048   &6  &0.331 &-9.416  &3.9387e-5 &11 &0.588 &-9.433 &3.9387e-5 \\
         &1/256    &261,120   &6  &0.666 &-9.216  &9.8111e-6 &11 &1.187 &-9.368 &9.8118e-6 \\
         &1/512    &523,264   &7  &1.736 &-9.936  &2.4178e-6 &11 &2.669 &-9.317 &2.4178e-6 \\
         &1/1024   &1,047,552 &7  &2.850 &-9.941  &7.4549e-7 &11 &4.396 &-9.300 &7.4557e-7 \\
\hline
\end{tabular}
\end{center}
\end{table}
According to Tables \ref{table:nonlin}--\ref{tab3},
we note that the used method with preconditioner $\mathcal{P}_{\alpha}$ converges much faster
than with $\mathcal{P}_1$ in terms of both \texttt{CPU} and \texttt{Iter} on this
Example~1, with different values of $\gamma$'s. In terms of accuracy (the values
\texttt{TRR} and \texttt{Err}), the two preconditioned Ad-SGMRES methods are
almost comparable. The results indicate that introducing the adaptive parameter $\alpha
\in(0,1)$ indeed helps improve the performance of $\mathcal{P}_1$. Moreover,
the $\tau$-matrix approximation of the Jacobian matrix $A$ in~(\ref{eq1.7x}) is
numerically effective. The iteration number of Ad-SGMRES-$\mathcal{P}_{\alpha}$
varies only slightly when $N_t$ increases (or $h$ decreases), showing almost
matrix-size independent convergence rate. Such favourable \textit{numerical scalability}
property of Ad-SGMRES-$\mathcal{P}_{\alpha}$ is completely in line with our theoretical analysis presented in Section~\ref{sec3.1}, where the eigenvalues distribution of the preconditioned matrices (also partly shown in Fig.~\ref{fig4.1}) is independent of the space-discretization matrix $A$. On the other hand, if the diffusion coefficient $\kappa_{\gamma}$ becomes smaller than $10^{-2}$, the difference of performance between the preconditioners
$\mathcal{P}_{\alpha}$ and $\mathcal{P}_1$ will be more remarkable; refer to Remark~\ref{rem3.2},
however we do not investigate these circumstance in the present study.

\vspace{1mm}
\noindent\textbf{Example 2.}~(the 2D problem \cite{Hao2015,Xu2019})
Consider the following Riesz fractional diffusion equation
\begin{equation}
\begin{cases}
\frac{\partial u(x,y,t)}{\partial t} = \kappa_{\gamma_1}\frac{\partial^{\gamma_1} u(x,y,t)
}{\partial |x|^{\gamma_1}} + \kappa_{\gamma_2}\frac{\partial^{\gamma_2u(x,y,t)}}{\partial
|y|^{\gamma_2}} + f(x,y,t), & (x,y,t)\in\Omega\times(0,T],\\
u(x,y,0) = \phi(x,y), & (x,y)\in\Omega,\\
u(x,y,t) = 0, & (x,y,t)\in\partial\Omega\times(0,T),
\label{eq4.1}
\end{cases}
\end{equation}
where $\kappa_{\gamma_1} = \kappa_{\gamma_2} = 0.01$, $T = 2$, $\Omega = (0,2)\times(0,2)$,
$\phi(x,y) = x^4(2 - x)^{4}y^4(2 - y)^4$ such that the source term
is exactly defined as
\begin{equation*}
\begin{split}
f(x,y,t) & = -\frac{1}{3}e^{-t/3}x^4(2 - x)^4y^4(2 - y)^4 + \frac{\kappa_{\gamma_1}e^{-t/3}}
{2\cos(\gamma_1\pi/2)}y^4(2 - y)^4\sum^{9}_{\ell = 5}\frac{q_{\ell}\Gamma(\ell)\big[x^{\ell
- 1 - \gamma_1} + (2 - x)^{\ell - 1 - \gamma_1}\big]}{\Gamma(\ell - \gamma_1)}\\
&\quad + \frac{\kappa_{\gamma_2}e^{-t/3}}{2\cos(\gamma_2\pi/2)}x^4(2 - x)^4\sum^{9}_{\ell =
5}\frac{q_{\ell}\Gamma(\ell)\big[x^{\ell - 1 - \gamma_2} + (2 - x)^{\ell - 1 - \gamma_2}\big]}
{\Gamma(\ell - \gamma_2)}
\end{split}
\end{equation*}
with $q_5 = 16,~q_6 = -32,~q_7 = 24,~q_8 = -8$ and $q_9 = 1$. The exact solution is
$u(x,y,y) = e^{-t/3}x^4(2 - x)^4y^4(2 - y)^4$.

By a similar derivation to the technique described in Section~\ref{sec2.1}, we can establish an implicit difference scheme for solving this model problem. For simplicity, we set $h_x = h_y = (b - a)/N$; then, it is easy to derive the Jacobian matrix in the following Kronecker product form
\begin{equation}
A = -\left[\left(\frac{\kappa_{\gamma_1}}{h^{\gamma_1}}T_x\right)\otimes I_{N-1} +
I_{N-1}\otimes \left(\frac{\kappa_{\gamma_2}}{h^{\gamma_2}}T_y\right)\right],
%
\label{eq4.2}
\end{equation}
where the two SPD Toeplitz matrices are defined by
\begin{equation*}
T_x = \begin{bmatrix}
\omega^{(\gamma_1)}_0 &\omega^{(\gamma_1)}_{-1} &\omega^{(\gamma_1)}_{-2} &\cdots
&\omega^{(\gamma_1)}_{3-N}&\omega^{(\gamma_1)}_{2-N}\\
\omega^{(\gamma_1)}_1 &\omega^{(\gamma_1)}_0&\omega^{(\gamma_1)}_{-1}&\cdots
&\omega^{(\gamma_1)}_{4-N}&\omega^{(\gamma_1)}_{3-N}\\
\omega^{(\gamma_1)}_2 &\omega^{(\gamma_1)}_1&\omega^{(\gamma_1)}_{0}&\cdots
&\omega^{(\gamma_1)}_{5-N}&\omega^{(\gamma_1)}_{4-N}\\
\vdots&\vdots&\vdots&\ddots&\vdots&\vdots\\
\omega^{(\gamma_1)}_{N-3} &\omega^{(\gamma_1)}_{N-4}&\omega^{(\gamma_1)}_{N-5}&\cdots
&\omega^{(\gamma_1)}_0&\omega^{(\gamma_1)}_{-1}\\
\omega^{(\gamma_1)}_{N-2} &\omega^{(\gamma_1)}_{N-3}&\omega^{(\gamma_1)}_{N-4}&\cdots
&\omega^{(\gamma_1)}_1&\omega^{(\gamma_1)}_0
\end{bmatrix}~~ {\rm and}~~
T_y = \begin{bmatrix}
\omega^{(\gamma_2)}_0 &\omega^{(\gamma_2)}_{-1} &\omega^{(\gamma_2)}_{-2} &\cdots
&\omega^{(\gamma_2)}_{3-N}&\omega^{(\gamma_2)}_{2-N}\\
\omega^{(\gamma_2)}_1 &\omega^{(\gamma_2)}_0&\omega^{(\gamma_2)}_{-1}&\cdots
&\omega^{(\gamma_1)}_{4-N}&\omega^{(\gamma_2)}_{3-N}\\
\omega^{(\gamma_2)}_2 &\omega^{(\gamma_2)}_1&\omega^{(\gamma_2)}_{0}&\cdots
&\omega^{(\gamma_2)}_{5-N}&\omega^{(\gamma_2)}_{4-N}\\
\vdots&\vdots&\vdots&\ddots&\vdots&\vdots\\
\omega^{(\gamma_2)}_{N-3} &\omega^{(\gamma_2)}_{N-4}&\omega^{(\gamma_2)}_{N-5}&\cdots
&\omega^{(\gamma_2)}_0&\omega^{(\gamma_2)}_{-1}\\
\omega^{(\gamma_2)}_{N-2} &\omega^{(\gamma_2)}_{N-3}&\omega^{(\gamma_2)}_{N-4}&\cdots
&\omega^{(\gamma_2)}_1&\omega^{(\gamma_2)}_0
\end{bmatrix}.
\end{equation*}
Then, it is not hard to know that $A$ is a SPD matrix, which has no impact on the algorithm
and theoretical analyses described in Sections \ref{sec2}--\ref{sec3}. The only difference
is that the direct application of $P_{\alpha}$ for Krylov subspace solvers must be prohibited,
because it requires the huge computational cost and memory storage for solving a sequence of dense shifted
linear systems in Step-(b).

In order to reduce the computational expense, we can follow Eq.~(\ref{eq3.9}) to approximate the Jacobian matrix $A$ in~(\ref{eq4.2}) by the following $\tau$-matrix
\begin{equation}
\begin{split}
\mathscr{T}(A) &= -\left[\left(\frac{\kappa_{\gamma_1}}{h^{\gamma_1}}\mathscr{T}(T_x)\right)\otimes I_{N-1} +
I_{N-1}\otimes \left(\frac{\kappa_{\gamma_2}}{h^{\gamma_2}}\mathscr{T}(T_y)\right)\right]\\
& = (Q^{\top}_{N-1}\otimes Q_{N-1})\left[\left(\frac{\kappa_{\gamma_1}}{h^{\gamma_1}}\Lambda^{(x)}_{N-1}\right)\otimes I_{N-1} +
I_{N-1}\otimes \left(\frac{\kappa_{\gamma_2}}{h^{\gamma_2}}\Lambda^{(y)}_{N-1}\right)\right](Q^{\top}_{N-1}\otimes Q_{N-1})
\end{split}
\end{equation}
where $\mathscr{T}(T_x) = Q^{\top}_{N-1}\Lambda^{(x)}_{N-1} Q_{N-1}$ and $\mathscr{T}(T_y) = (Q^{\top}_{N-1}\Lambda^{(y)}_{N-1} Q_{N-1})$ \cite{Bini90,Nout11,Serra99}. Again, instead of solving the shifted linear systems in Step-(b), we solve the following sequence of shifted linear systems,
\begin{eqnarray}
\big[\lambda^{(\alpha)}_{n} I_s - \tau\mathscr{T}(A)\big]{\bm z}_{2,n} = {\bm z}_{1,n},
\quad n = 1,2,\cdots,N_t,\\
\Leftrightarrow(Q^{\top}_{N-1}\otimes Q_{N-1})\left[\lambda^{(\alpha)}_n
I_{s} - \left(\frac{\kappa_{\gamma_1}}{h^{\gamma_1}}\Lambda^{(x)}_{N-1}\right)\otimes I_{N-1} +
I_{N-1}\otimes \left(\frac{\kappa_{\gamma_2}}{h^{\gamma_2}}\Lambda^{(y)}_{N-1}\right)\right]
(Q^{\top}_{N-1}\otimes Q_{N-1}){\bm z}_{2,n} = {\bm z}_{1,n},\label{eq4.5}
\end{eqnarray}
where $s = (N-1)^{2}$, efficiently by fast discrete sine transforms avoiding the storage of any dense matrices. The numerical results reported in Tables~\ref{tab4}--\ref{tab6} with $h = h_x = h_y$ illustrate the effectiveness and robustness of the PinT preconditioner $\mathcal{P}_{\alpha}$ for solving~Eq.~(\ref{eq4.5}).
%
\begin{table}[!htpb]\tabcolsep=5pt
\begin{center}
\caption{Numerical results of GMRES with two different preconditioners on Example 2
with $(\gamma_1,\gamma_2) = (1.4,1.2)$.}
\vspace{1mm}
\begin{tabular}{crrccrccccc}
\hline &&&\multicolumn{4}{c}{$\mathcal{P}_{\alpha}$} &\multicolumn{4}{c}{$\mathcal{P}_1$}\\
[-2pt]\cmidrule(l{0.7em}r{0.7em}){4-7} \cmidrule(l{0.7em}r{0.6em}){8-11}\\[-11pt]
$N_t$ &$h_x = h_y$ &DoF &\texttt{Iter} &\texttt{CPU} &\texttt{TRR} &\texttt{Err}
&\texttt{Iter}&\texttt{CPU} &\texttt{TRR}&\texttt{Err} \\
\hline
$2^{6}$  &1/64  &254,016     &4.0 &1.006  &-9.195  &1.2627e-4 &12.0 &2.544  &-9.290 &1.2628e-4 \\
         &1/128 &1,032,256   &4.5 &5.325  &-10.066 &8.0645e-5 &12.0 &12.312 &-9.132 &8.0646e-5 \\
         &1/256 &4,161,600   &4.5 &17.597 &-9.219  &7.8998e-5 &12.0 &42.376 &-9.134 &7.8999e-5 \\
         &1/512 &16,711,744  &5.0 &115.94 &-9.814  &7.8611e-5 &12.5 &262.03 &-9.151 &7.8612e-5 \\
\hline
$2^8$    &1/64  &1,016,064   &4.0 &4.027  &-9.544  &7.4633e-5 &12.0 &10.355 &-9.264 &7.4633e-5 \\
         &1/128 &4,129,024   &4.0 &19.419 &-9.134  &2.1246e-5 &12.0 &49.827 &-9.102 &2.1246e-5 \\
         &1/256 &16,646,400  &5.0 &79.796 &-9.446  &7.8953e-6 &12.0 &173.68 &-9.092 &7.8954e-6 \\
         &1/512 &66,846,976  &4.5 &426.03 &-9.174  &5.0729e-6 &12.5 &1050.3 &-9.204 &5.0732e-6 \\
\hline
$2^{10}$ &1/64  &4,064,256   &4.0 &16.246 &-10.084 &7.1404e-5 &12.0 &41.638 &-9.244 &7.1404e-5 \\
         &1/128 &16,516,096  &4.0 &78.194 &-9.861  &1.8016e-5 &12.0 &203.61 &-9.084 &1.8017e-5 \\
         &1/256 &66,585,600  &4.0 &265.24 &-9.567  &4.6657e-6 &12.0 &753.81 &-9.072 &4.6661e-6 \\
         &1/512 &267,387,904 &4.0 &2520.7 &-9.157  &1.3276e-6 &12.5 &6989.3 &-9.209 &1.3282e-6 \\
\hline
\end{tabular}
\end{center}
\label{tab4}
\end{table}
\begin{table}[t]
\begin{center}
\caption{Numerical results of GMRES with two different preconditioners on Example 2
with $(\gamma_1,\gamma_2) = (1.5,1.5)$.}
\vspace{1mm}
\begin{tabular}{crrcccccccc}
\hline &&&\multicolumn{4}{c}{$\mathcal{P}_{\alpha}$} &\multicolumn{4}{c}{$\mathcal{P}_1$}\\
[-2pt]\cmidrule(l{0.7em}r{0.7em}){4-7} \cmidrule(l{0.7em}r{0.6em}){8-11}\\[-11pt]
$N_t$ &$h_x = h_y$ &DoF &\texttt{Iter} &\texttt{CPU} &\texttt{TRR} &\texttt{Err}
&\texttt{Iter} &\texttt{CPU} &\texttt{TRR} &\texttt{Err} \\
\hline
$2^{6}$  &1/64  &254,016     &4.0 &1.004  &-9.297  &1.5758e-4 &11.0 &2.487  &-9.255 &1.5758e-4 \\
         &1/128 &1,032,256   &4.5 &5.322  &-9.978  &8.1963e-5 &11.0 &11.416 &-9.227 &8.1963e-5 \\
         &1/256 &4,161,600   &5.0 &19.756 &-10.122 &7.9075e-5 &11.5 &41.569 &-9.379 &7.9075e-5 \\
         &1/512 &16,711,744  &5.0 &115.82 &-9.676  &7.8490e-5 &11.5 &240.18 &-9.345 &7.8490e-5 \\
\hline
$2^8$    &1/64  &1,016,064   &4.0 &4.058  &-9.585  &1.0778e-4 &11.0 &9.745  &-9.127 &1.0778e-4 \\
         &1/128 &4,129,024   &4.0 &19.609 &-9.106  &2.9441e-5 &11.0 &46.275 &-9.162 &2.9441e-5 \\
         &1/256 &16,646,400  &4.5 &73.342 &-9.454  &9.8515e-6 &11.5 &168.30 &-9.259 &9.8515e-6 \\
         &1/512 &66,846,976  &4.5 &424.74 &-9.300  &5.1188e-6 &11.5 &965.63 &-9.245 &5.1183e-6 \\
\hline
$2^{10}$ &1/64  &4,064,256   &4.0 &16.343 &-10.060 &1.0466e-4 &11.0 &39.182 &-9.100 &1.0466e-4 \\
         &1/128 &16,516,096  &4.0 &78.225 &-9.813  &2.6328e-5 &11.0 &188.46 &-9.126 &2.6328e-5 \\
         &1/256 &66,585,600  &4.0 &265.75 &-9.582  &6.7380e-6 &11.5 &733.89 &-9.256 &6.7382e-6 \\
         &1/512 &267,387,904 &4.0 &2517.9 &-9.282  &1.8400e-6 &11.5 &6278.1 &-9.224 &1.8401e-6 \\
\hline
\end{tabular}
\end{center}
\label{tab5}
\end{table}
\begin{table}[!htpb]
\begin{center}
\caption{Numerical results of GMRES with two different preconditioners on Example 2
with $(\gamma_1,\gamma_2) = (1.7,1.9)$.}
\label{tab6}
\vspace{1mm}
\begin{tabular}{crrcccccccc}
\hline &&&\multicolumn{4}{c}{$\mathcal{P}_{\alpha}$} &\multicolumn{4}{c}{$\mathcal{P}_1$}\\
[-2pt]\cmidrule(l{0.7em}r{0.7em}){4-7} \cmidrule(l{0.7em}r{0.6em}){8-11}\\[-11pt]
$N_t$ &$h_x = h_y$ &DoF &\texttt{Iter} &\texttt{CPU} &\texttt{TRR} &\texttt{Err}
&\texttt{Iter} &\texttt{CPU} &\texttt{TRR} &\texttt{Err}\\
\hline
$2^{6}$  &1/64   &254,016     &4.0 &1.047  &-9.492  &2.3321e-4 &11.5 &2.513  &-10.177 &2.3321e-4\\
         &1/128  &1,032,256   &4.0 &4.976  &-9.058  &9.5250e-5 &11.5 &12.175 &-10.188 &9.5251e-5 \\
         &1/256  &4,161,600   &4.5 &17.826 &-9.887  &7.9355e-5 &11.5 &41.627 &-10.208 &7.9355e-5 \\
         &1/512  &16,711,744  &4.5 &114.13 &-9.473  &7.8240e-5 &11.0 &231.76 &-9.154  &7.8240e-5 \\
\hline
$2^8$    &1/64   &1,016,064   &4.0 &4.035  &-10.062 &1.8715e-4 &11.0 &9.673  &-9.339  &1.8715e-4 \\
         &1/128  &4,129,024   &4.0 &19.567 &-9.425  &4.9102e-5 &11.5 &48.305 &-9.977  &4.9102e-5 \\
         &1/256  &16,646,400  &4.0 &66.497 &-9.150  &1.4579e-5 &11.0 &161.38 &-9.064  &1.4578e-5 \\
         &1/1024 &66,846,976  &4.0 &380.78 &-9.011  &5.9522e-6 &11.0 &921.53 &-9.115  &5.9518e-6 \\
\hline
$2^{10}$ &1/64   &4,064,256   &4.0 &16.225 &-10.092 &1.8428e-4 &11.0 &38.684 &-9.027  &1.8428e-4 \\
         &1/128  &16,516,096  &4.0 &78.158 &-9.905  &4.6224e-5 &11.5 &198.25 &-10.077 &4.6224e-5 \\
         &1/256  &66,585,600  &4.0 &266.02 &-9.760  &1.1700e-5 &11.5 &734.78 &-10.020 &1.1701e-5 \\
         &1/512  &267,387,904 &4.0 &2519.5 &-9.551  &3.0690e-6 &11.0 &6117.6 &-9.100  &3.0691e-6 \\
\hline
\end{tabular}
\end{center}
\end{table}

Similar to Example 1, Tables \ref{tab4}--\ref{tab6} show that the preconditioner
$\mathcal{P}_{\alpha}$ converges much faster in terms of both \texttt{CPU} and
\texttt{Iter} than $\mathcal{P}_1$ on this example, with different values of $
\gamma$'s. The accuracy of two preconditioned BiCGSTAB methods is almost comparable
with respect to \texttt{TRR} and \texttt{Err}. Once again, the results indicate that
introducing the adaptive parameter $\alpha\in(0,1)$ indeed helps improve the performance of
$\mathcal{P}_1$; the faster convergence rate of $\mathcal{P}_{\alpha}$ is computational
attractive especially when the diffusion coefficients become smaller than $10^{-2}$ -- cf.
Remark \ref{rem3.2}. The $\tau$-matrix approximation to the Jacobian matrix $A$ in~(\ref{eq1.7x}) remains very effective for this two-dimensional model problem~(\ref{eq4.1}). For more details, one can attempt to play it by our accessible MATLAB codes. The iteration number of BiCGSTAB-$\mathcal{P}_{\alpha}$ varies only slightly when $N_t$ increases (or $h$ decreases),
showing also for this problem almost matrix-size independent convergence rate.
Such favourable \textit{numerical scalability} property confirms our theoretical
analysis since in Section~\ref{sec3.1} the eigenvalue distribution of the preconditioned
matrices is independent of the space-discretization matrix $A$.
\section{Conclusions}
\label{sec5}
In this note, we revisit the all-at-once linear system arising from the BDF$p$ temporal discretization for evolutionary PDEs. In particular, we present the BDF2 scheme for the
RFDEs model as our case study, where the resultant all-at-once system is a BLTT linear system with a low-rank matrix perturbation. The conditioning of the all-at-once coefficient matrix is studied and tight bounds are provided. Then, we adapt the generalized BC preconditioner for
such all-at-once systems,  proving the invertibility of the preconditioner matrix unlike in previous studies. Our analysis demonstrates the superiority of the generalized BC preconditioner to the BC preconditioner. Moreover, the spectral properties of the preconditioned system and the convergence behavior of the preconditioned Krylov subspace solver have been investigated. By the $\tau$-matrix approximation of the dense Jacobian matrix $A$, we derive a memory-effective implementation of the generalized BC preconditioner for solving the one- and two-dimensional RFDEs model problems. Numerical results have been reported to show the effectiveness of the generalized BC preconditioner.

On the other hand, according to the analysis and implementation of the generalized BC
preconditioner, the BDF2 scheme implemented via PinT preconditioned Krylov solvers
can be easily extended to other spatial discretizations schemes for RFDEs (with other boundary conditions). Furthermore, our study may inspire the development of new parallel
numerical methods preserving the positivity for certain evolutionary PDEs. As an outlook for the future, the extension of the generalized BC preconditioner and of its parallel implementation for solving nonlinear RFDEs \cite{Gander17,Zhao19} (even in cases when non-uniform temporal steps are chosen in combination with the BDF2 scheme; refer to \cite{Goddard,Falgout} for a short discussion) remains an interesting topic of further research.

\section*{Acknowledgments}
{\em The authors would like to thank Dr. Xin Liang and Prof. Shu-Lin
Wu for their insightful suggestions and encouragements. This research
is supported by NSFC (11801463 and 61876203) and the Applied Basic Research Project
of Sichuan Province (2020YJ0007). Meanwhile, the first author is
grateful to Prof. Hai-Wei Sun for his helpful discussions during the
visiting to University of Macau.}
\vspace{5mm}


\end{document}